\documentclass[10pt,reqno]{amsart}

\usepackage{amsmath, epsfig, cite}
\usepackage{amssymb}
\usepackage{amsfonts}
\usepackage{latexsym}
\usepackage{graphicx}

\newtheorem{thm}{Theorem}[section]

\def\pf{\noindent {\it Proof.} }

\numberwithin{equation}{section}

\makeatletter \@addtoreset{equation}{section} \makeatother

\begin{document}
\allowdisplaybreaks

\title
{One-Parameter Generalizations of Rogers-Ramanujan Type Identities}
\author[N. S. S. Gu]{Nancy S. S. Gu\textsuperscript{\dag}}
\address[N. S. S. Gu]{Center for Combinatorics, LPMC-TJKLC, Nankai
University, Tianjin 300071, P. R. China} \email{gu@nankai.edu.cn}

\author[H. Prodinger]{Helmut Prodinger\textsuperscript{*}}
\address[H. Prodinger]{Department of Mathematics\\
University of Stellenbosch\\
7602 Stellenbosch\\
South Africa} \email{hproding@sun.ac.za}
\thanks{\textsuperscript{\dag} The first author was supported by the PCSIRT
Project of the Ministry of Education, and the Specialized Research
Fund for the Doctoral Program of Higher Education of China
(200800551042). }
\thanks{\textsuperscript{*} The second author was supported by International
Science and Technology Agreement (Grant 67215).}

\keywords{Rogers-Ramanujan identities, $q$-series, determinant,
recursion.} \subjclass{05A30, 15A15}
\date{\today}
\begin{abstract}
Resorting to the recursions satisfied by the polynomials which
converge to the right hand sides of the Rogers-Ramanujan type
identities given by Sills \cite{Sills} and a determinant method
presented in \cite{Ismail-Prodinger-Stanton}, we obtain many new
one-parameter generalizations of the Rogers-Ramanujan type
identities, such as a generalization of the analytic versions of the
first and second G\"{o}llnitz-Gordon partition identities, and
generalizations of the first, second, and third Rogers-Selberg
identities.
\end{abstract}

\maketitle

%==============================================================================

\section{Introduction} 

In \cite{Garrett-Ismail-Stanton}, by evaluating an integral
involving $q$-Hermite polynomials in two different ways and equating
the results, Garrett et al. found a generalization of the celebrated
Rogers-Ramanujan identities:
\begin{equation}\label{GIS}
\sum_{n=0}^{\infty}\frac{q^{n^2+mn}}{(q;q)_n}=\frac{(-1)^mq^{-{m
\choose 2}}E_{m-2}}{(q,q^4;q^5)_{\infty}}-\frac{(-1)^mq^{-{m \choose
2}}D_{m-2}}{(q^2,q^3;q^5)_{\infty}},
\end{equation}
where the Schur polynomials $D_m$ and $E_m$ are defined by
\begin{align*}
&D_m=D_{m-1}+q^mD_{m-2}, \qquad D_0=1, \  D_1=1+q, \\
&E_m=E_{m-1}+q^mE_{m-2}, \qquad E_0=1, \  E_1=1,
\end{align*}
and Schur \cite{Schur} gave the limit
\begin{equation*}
D_{\infty}=\frac{1}{(q,q^4;q^5)_{\infty}}, \qquad
E_{\infty}=\frac{1}{(q^2,q^3;q^5)_{\infty}}.
\end{equation*}

It is obvious that we can get the following two Rogers-Ramanujan
identities by letting $m=0$ and $m=1$ in \eqref{GIS}, respectively.
\begin{equation}\label{Rogers1}
\sum_{n=0}^{\infty}\frac{q^{n^2}}{(q;q)_n}=\frac{1}{(q,q^4;q^5)_{\infty}},
\end{equation}
\begin{equation}\label{Rogers2}
\sum_{n=0}^{\infty}\frac{q^{n^2+n}}{(q;q)_n}=\frac{1}{(q^2,q^3;q^5)_{\infty}}.
\end{equation}

Later, Andrews et al. \cite{Andrews-Knopfmacher-Paule} provided an
alternative proof of \eqref{GIS} by using the extended Engel
expansion. In \cite{Ismail-Prodinger-Stanton}, Ismail et al. used
the theory of associated orthogonal polynomials to explain
determinants that Schur introduced in 1917, and showed that Equation
\eqref{GIS} can be obtained from the Rogers-Ramanujan identities
\eqref{Rogers1} and \eqref{Rogers2}. Furthermore, Andrews et al.
\cite{Andrews-Knopfmacher-Paule-Prodinger} discussed Al-Salam/Ismail
and Santos polynomials in the context of identities of \eqref{GIS}
type.

The main purpose of this paper is to apply the determinant method
which was presented in \cite{Ismail-Prodinger-Stanton} to generalize
the Rogers-Ramanujan type identities. In \cite{Sills}, Sills mainly
focused on a method which was developed by Andrews \cite[$\S$9.2, p.
88]{Andrews} for discovering finite analogs of Rogers-Ramanujan type
identities via $q$-difference equations. In the paper, he presented
at least one finitization for each of the $130$ identities in
Slater's list \cite{Slater}, along with recursions satisfied by the
polynomials which converge to the right hand sides of the
Rogers-Ramanujan type identities. Resorting to these recursions and
the determinant method, we obtain many new parameterized
generalizations of the Rogers-Ramanujan type identities, such as a
generalization of the analytic versions of the first and second
G\"{o}llnitz-Gordon partition identities, and generalizations of the
first, second, and third Rogers-Selberg identities. In Section $2$,
we mainly discuss the three-term recursions. In Section $3$, we
focus on four-term recursions. Moreover, in
\cite{Bowman-McLaughin-Sills,McLaughlin-Sills-Zimmer}, the authors
also found some new Rogers-Ramanujan type identities which are the
partners to those in Slater's list. By using the determinant method,
we can give the initial conditions of the recursions for these new
identities, and then find the generalizations of these identities.

In \cite{Sills}, Sills gave an annotated and cross-referenced
version of Slater's list of identities from \cite{Slater} as an
appendix. In this paper, we use this version of the list as the
reference.

As usual, we follow the notation and terminology in
\cite{Gasper-Rahman}. For $|q|<1$, the $q$-shifted factorial is
defined by
$$(a;q)_\infty=
\prod_{k=0}^{\infty}(1-aq^k) \text{\ \  and \ \  }(a;q)_n
=\frac{(a;q)_\infty}{(aq^n;q)_\infty}, \text{ for } n\in
\mathbb{C}.$$

For convenience, we shall adopt the following notation for multiple
$q$-shifted factorials:
$$(a_1,a_2,\ldots,a_m;q)_n=(a_1;q)_n(a_2;q)_n\cdots(a_m;q)_n,$$
where $n$ is an integer or infinity.

In order to sketch the paper clearly, we list the main results in a
table.

\newpage

\begin{center}
\begin{tabular}{|l|c|}
\hline
\textbf{Identities in Slater's list and some new ones} & \textbf{Generalizations} \\
\hline

Identity A.8 (Gauss-Lebesgue \cite{Lebesgue}) & Theorem
\ref{A8-A13Thm} \\[3pt]
Identity A.13 (Slater \cite{Slater}) &  \\[3pt]

\hline

Identity A.16 (Rogers \cite{Rogers}) & Theorem
\ref{A16-A20Thm}\\[3pt]
Identity A.20 (Rogers \cite{Rogers})& \\[3pt]

\hline

Identity A.29 (Slater \cite{Slater}) & Theorem
\ref{A29-A50Thm} \\[3pt]
Identity A.50 (Slater \cite{Slater}) & (1) (2) \\[3pt]

\hline

Identity A.34 (Slater \cite{Slater}): The analytic version of the
second &  \\
G\"{o}llnitz-Gordon partition identity. & Theorem
\ref{A34-A36Thm}\\[3pt]
Identity A.36 (Slater \cite{Slater}): The analytic version of the first & \\
G\"{o}llnitz-Gordon partition identity. & \\[3pt]

\hline

Identity A.38 (Slater \cite{Slater}) & Theorem \ref{A38-A39Thm} \\[3pt]
Identity A.39 (Jackson \cite{Jackson}) & (1) (2) \\[3pt]

\hline

Identity A.79 (Rogers \cite{Rogers}) & Theorem \ref{A79-A96Thm} \\[3pt]
Identity A.96 (Rogers \cite{Rogers}) & (1) (2)\\[3pt]

\hline

Identity A.94 (Rogers \cite{Rogers}) & Theorem \ref{A94-A99Thm} \\[3pt]
Identity A.99 (Rogers \cite{Rogers}) & (1) (2)\\[3pt]

\hline

Identity A.25 (Slater \cite{Slater}) & Theorem \ref{A25-A25'Thm} \\[3pt]
An identity (McLaughlin et al. \cite[Eq.
(2.7)]{McLaughlin-Sills-Zimmer})& \\[3pt]

\hline

Identity A.31 (Rogers \cite{Rogers1} and Selberg
\cite{Selberg}) & \\
The third Rogers-Selberg identity &  \\[3pt]
Identity A.32 (Rogers \cite{Rogers} and Selberg
\cite{Selberg}) & Theorem \ref{A31-A32-A33Thm}\\
The second Rogers-Selberg identity &  (1) (2)\\[3pt]
Identity A.33 (Rogers \cite{Rogers} and Selberg \cite{Selberg})
&\\
The first Rogers-Selberg identity & \\[3pt]

\hline

Identity A.59 (Rogers \cite{Rogers1}) & \\[3pt]

Identity A.60 (Rogers \cite{Rogers1}) & Theorem \ref{A59-A60-A61Thm}\\[3pt]

Identity A.61 (Rogers \cite{Rogers}) & (1) (2)\\[3pt]

\hline

Identity A.80 (Rogers \cite{Rogers1}) & \\[3pt]

Identity A.81 (Rogers \cite{Rogers1}) & Theorem \ref{A80-A81-A82Thm}\\[3pt]

Identity A.82 (Rogers \cite{Rogers1}) & (1) (2)\\[3pt]

\hline

Identity A.117 (Slater \cite{Slater}) & \\[3pt]

Identity A.118 (Slater \cite{Slater}) & Theorem \ref{A117-A118-A119Thm}\\[3pt]

Identity A.119 (Slater \cite{Slater}) & (1) (2)\\[3pt]

\hline

Identity A.21 (Slater \cite{Slater}) & \\[3pt]

An identity (McLaughlin et al. \cite[Eq.
(2.5)]{McLaughlin-Sills-Zimmer}) & Theorem \ref{A21-Zimmer2.5-Bowman2.7Thm}\\[3pt]

An identity (Bowman et al. \cite[Thm. 2.7]{Bowman-McLaughin-Sills}) & (1) (2)\\[3pt]

\hline
\end{tabular}
\end{center}

\section{Generalizations of identities with three-term recursions}

In this section, we generalize the Rogers-Ramanujan type identities
in Slater's list \cite{Slater} by using the determinant method
presented in \cite{Ismail-Prodinger-Stanton}. Start with the
three-term recursions of the polynomials which converge to the right
hand sides of the identities in \cite{Sills}. First, we construct a
function $F(z)$ which is expressed by an infinite determinant. By
expanding the determinant and comparing the coefficients, we get a
summation expression of $F(z)$. Then, we expand $D_n(z)$, a finite
determinant of $F(z)$, to get a recursion which has appeared in
Sills' list \cite[Sec. 3.2]{Sills}. Assume that the polynomials
$P_n$ and $Q_n$ satisfy this recursion with different initial
conditions, then $D_n(z)$ can be expressed by a linear combination
of these two polynomials. By means of the initial conditions of
$D_n(z)$, we get the limit of $D_n(z)$ which is another expression
of $F(z)$. Finally, equating the two different expressions of
$F(z)$, we obtain a new generalization.

In the following, for convenience, the recursions given by Sills
\cite{Sills} are directly presented below the identities in Slater's
list.

\begin{thm}\label{A8-A13Thm}We have
\begin{equation}\label{A8-A13}
\sum_{n=0}^{\infty}\frac{(-q;q)_n
q^{n(n+2m-1)/2}}{(q;q)_n}=(-1)^mq^{-{m \choose
2}}Q_{m-1}\frac{(-q;q^2)_{\infty}}{(q;q^2)_{\infty}}-(-1)^mq^{-{m
\choose 2}}R_{m-1}\frac{(q^4;q^4)_{\infty}}{(q;q)_{\infty}},
\end{equation}
where
\begin{align*}
Q_m &=(1+q^{m-1})Q_{m-1}+q^{m-1}Q_{m-2},\qquad Q_{-1}=1, \ Q_0=0, \
Q_1=1, \\
R_m &=(1+q^{m-1})R_{m-1}+q^{m-1}R_{m-2},\qquad R_{-1}=-1,\ R_0=1, \
R_1=1.
\end{align*}
\end{thm}

\pf The identities A.8 and A.13 in Slater's list are stated as
follows.

\textbf{Identity A.8 (Gauss-Lebesgue \cite{Lebesgue}):}

\begin{equation}\label{A8}
\sum_{n=0}^{\infty}\frac{(-q;q)_n q^{n(n+1)/2}}{(q;q)_n}=
\frac{(q^4;q^4)_{\infty}}{(q;q)_{\infty}}.
\end{equation}
Sills \cite{Sills} gave the following recursion for $P_n$ which
converge to the right hand side of \eqref{A8}.
\begin{equation}\label{A8P}
P_n=(1+q^n)P_{n-1}+q^nP_{n-2},\qquad P_{-1}=0,\ P_0=1, \ P_1=1+q.
\end{equation}

\textbf{Identity A.13 (Slater \cite{Slater}):}

\begin{equation}\label{A13}
\sum_{n=0}^{\infty}\frac{(-q;q)_n q^{n(n-1)/2}}{(q;q)_n}=
\frac{(q^4;q^4)_{\infty}}{(q;q)_{\infty}}+\frac{(-q;q^2)_{\infty}}{(q;q^2)_{\infty}},
\end{equation}
\begin{equation}\label{A13P}
P_n=(1+q^{n-1})P_{n-1}+q^{n-1}P_{n-2},\qquad P_{-1}=0,\ P_0=1, \
P_1=2.
\end{equation}

First, we need to shift the index $n$ in \eqref{A8P} to let the two
recursions coincide with each other. Letting $Q_n=P_{n-1}$ in
\eqref{A8P}, we get
\begin{equation}\label{A8Q}
Q_n=(1+q^{n-1})Q_{n-1}+q^{n-1}Q_{n-2},\qquad Q_{-1}=1,\ Q_0=0, \
Q_1=1.
\end{equation}
Thus, $P_n$ in \eqref{A13P} and $Q_n$ in \eqref{A8Q} satisfy the
same recursion with different initial conditions, and converge to
the right hand sides of \eqref{A13} and \eqref{A8}, respectively. In
the following, we use $P_n$ in \eqref{A13P} and $Q_n$ in \eqref{A8Q}
to prove this theorem.

Then consider the following determinant:
\begin{equation*}
F(z):=\left| \begin{array}{ccccc} 1+z & zq & & & \cdots \\
-1 & 1+zq & zq^2 & & \cdots  \\
& -1 & 1+zq^2 & zq^3 & \cdots \\
& & \ddots & \ddots & \ddots\\
\end{array}
\right|.
\end{equation*}
Expanding the determinant with respect to the first column, we get
\begin{equation*}
F(z)=(1+z)F(zq)+zqF(zq^2).
\end{equation*}
Setting
\begin{equation*}
F(z)=\sum_{n=0}^{\infty}a_nz^n,
\end{equation*}
by comparing coefficients, we have
\begin{equation*}
a_n=q^na_n+q^{n-1}a_{n-1}+q^{2n-1}a_{n-1},
\end{equation*}
\begin{equation*}
a_n=\frac{(1+q^n)q^{n-1}}{1-q^n}a_{n-1}=\cdots=\frac{(-q;q)_nq^{n(n-1)/2}}{(q;q)_n}a_0.
\end{equation*}
Since $a_0=F(0)=1$, iteration leads to
\begin{equation*}
F(z)=\sum_{n=0}^{\infty}\frac{(-q;q)_nq^{n(n-1)/2}}{(q;q)_n}z^n,
\end{equation*}
and thus the left hand side of \eqref{A8-A13} can be expressed by
$F(q^m)$.

On the other hand, $F(z)$ is the limit of the finite determinant
\begin{align*}
D_n(z):= \left| \begin{array}{ccccc} 1+z & zq& & & \cdots \\
-1 & 1+zq & zq^2 & & \cdots \\
& -1 & 1+zq^2 & zq^3 & \cdots \\
\vdots & \vdots & \ddots & \ddots & \ddots \\
& & -1 & 1+zq^{n-2} & zq^{n-1} \\
& & & -1 & 1+zq^{n-1} \\
\end{array} \right|.
\end{align*}
Expanding this determinant with respect to the last row, we get
\begin{equation*}
D_n(z)=(1+zq^{n-1})D_{n-1}(z)+zq^{n-1}D_{n-2}(z), \qquad D_0(z)=1,\
D_1(z)=1+z.
\end{equation*}
Then we have
\begin{equation}\label{A8-A13D}
D_{n-m}(q^m)=(1+q^{n-1})D_{n-m-1}(q^m)+q^{n-1}D_{n-m-2}(q^m).
\end{equation}
According to \eqref{A13P}, \eqref{A8Q}, and \eqref{A8-A13D}, we
notice that the sequences $\langle D_{n-m}(q^m)\rangle_n$, $\langle
P_n\rangle_n$, and $\langle Q_n\rangle_n$ satisfy the same
recursion. Set
\begin{equation*}
D_{n-m}(q^m)=\lambda_m P_n+\mu_m Q_n.
\end{equation*}
We can determine the parameters $\lambda_m$ and $\mu_m$ using the
initial conditions $D_{0}(q^m)=1$, $D_{1}(q^m)=1+q^{m}$, and the
recursions \eqref{A13P} and \eqref{A8Q}, which leads to the
evaluations
\begin{equation*}
\lambda_m=\frac{Q_{m-1}}{P_mQ_{m-1}-P_{m-1}Q_m},
\end{equation*}
\begin{equation*}
\mu_m=\frac{P_{m-1}}{P_{m-1}Q_m-P_mQ_{m-1}}.
\end{equation*}
Indeed, we have
\begin{equation*}
P_mQ_{m-1}-P_{m-1}Q_m=(-1)^mq^{m \choose 2},
\end{equation*}
which can be proved by induction on $m$.

Therefore, we have simpler forms for $\lambda_m$ and $\mu_m$ as
follows:
\begin{equation*}
\lambda_m=(-1)^mq^{-{m \choose 2}}Q_{m-1},\hskip 10mm
\mu_m=-(-1)^mq^{-{m \choose 2}}P_{m-1}.
\end{equation*}
Notice that the above analysis has led to
\begin{equation*}
D_{n-m}(q^m)=(-1)^mq^{-{m \choose 2}}Q_{m-1}P_n-(-1)^mq^{-{m \choose
2}}P_{m-1}Q_n.
\end{equation*}
Letting $n\rightarrow \infty$, we have
\begin{equation*}
F(q^m)=(-1)^mq^{-{m \choose 2}}Q_{m-1}P_{\infty}-(-1)^mq^{-{m
\choose 2}}P_{m-1}Q_{\infty},
\end{equation*}
which is equivalent to the following identity
\begin{equation*}
\sum_{n=0}^{\infty}\frac{(-q;q)_n
q^{n(n+2m-1)/2}}{(q;q)_n}=(-1)^mq^{-{m \choose
2}}\frac{(-q;q^2)_{\infty}}{(q;q^2)_{\infty}}Q_{m-1}-(-1)^mq^{-{m
\choose
2}}\frac{(q^4;q^4)_{\infty}}{(q;q)_{\infty}}(P_{m-1}-Q_{m-1}).
\end{equation*}

Finally, set $R_{m-1}=P_{m-1}-Q_{m-1}$. According to \eqref{A13P}
and \eqref{A8Q}, we have
\begin{equation*}
R_m=(1+q^{m-1})R_{m-1}+q^{m-1}R_{m-2},\qquad R_{-1}=-1,\ R_0=1,\
R_1=1.
\end{equation*}
Therefore, we obtain \eqref{A8-A13} as desired.\qed

Setting $m=1$ and $m=0$ in \eqref{A8-A13}, we get the identities
\eqref{A8} and \eqref{A13}, respectively.

\begin{thm}\label{A16-A20Thm}We have
\begin{equation}\label{A16-A20}
\sum_{n=0}^{\infty}\frac{q^{n^2+2mn}}{(q^4;q^4)_n}=\frac{A_m}{(q,q^4;q^5)_{\infty}(-q^2;q^2)_{\infty}}
+\frac{B_m}{(q^2,q^3;q^5)_{\infty}(-q^2;q^2)_{\infty}},
\end{equation}
where
\begin{equation*}
A_m=-q^{2m-3}A_{m-1}+A_{m-2},\qquad A_0=1,\ A_1=0,
\end{equation*}
\begin{equation*}
B_m=-q^{2m-3}B_{m-1}+B_{m-2},\qquad B_0=0,\ B_1=1.
\end{equation*}
\end{thm}

\pf We state the identities A.16 and A.20 in Slater's list with the
recursions given by Sills \cite{Sills} as follows.

\textbf{Identity A.16 (Rogers \cite{Rogers}):}

\begin{equation}\label{A16}
\sum_{n=0}^{\infty}\frac{q^{n^2+2n}}{(q^4;q^4)_n}=
\frac{1}{(q^2,q^3;q^5)_{\infty}(-q^2;q^2)_{\infty}},
\end{equation}
\begin{equation}\label{A16P}
P_n=(1-q^2+q^{2n+1})P_{n-1}+q^2P_{n-2},\qquad P_{-1}=1,\ P_0=1, \
P_1=1+q^3.
\end{equation}

\textbf{Identity A.20 (Rogers \cite{Rogers}):}

\begin{equation}\label{A20}
\sum_{n=0}^{\infty}\frac{q^{n^2}}{(q^4;q^4)_n}=
\frac{1}{(q,q^4;q^5)_{\infty}(-q^2;q^2)_{\infty}},
\end{equation}
\begin{equation}\label{A20P}
P_n=(1-q^2+q^{2n-1})P_{n-1}+q^2P_{n-2},\qquad P_{-1}=1,\ P_0=1, \
P_1=1+q.
\end{equation}

For the recursion \eqref{A16P}, letting $Q_n=P_{n-1}$, we get
\begin{equation}\label{A16Q}
Q_n=(1-q^2+q^{2n-1})Q_{n-1}+q^2Q_{n-2},\qquad Q_{-1}=1-q^{-1},\
Q_0=1, \ Q_1=1.
\end{equation}

Therefore, $P_n$ in \eqref{A20P} and $Q_n$ in \eqref{A16Q} satisfy
the same recursion with different initial conditions and converge to
the right hand sides of \eqref{A20} and \eqref{A16}, respectively.

Consider the following determinant:
\begin{equation*}
F(z):=\left| \begin{array}{ccccc} 1-q^2+zq & q^2 & & & \cdots \\
-1 & 1-q^2+zq^3 & q^2 & & \cdots  \\
& -1 & 1-q^2+zq^5 & q^2 & \cdots \\
& & \ddots & \ddots & \ddots\\
\end{array}
\right|.
\end{equation*}
Expanding the determinant with respect to the first column, we get
\begin{equation*}
F(z)=(1-q^2+zq)F(zq^2)+q^2F(zq^4).
\end{equation*}
Setting
\begin{equation*}
F(z)=\sum_{n=0}^{\infty}a_nz^n,
\end{equation*}
we obtain, upon comparing coefficients,
\begin{equation*}
a_n=q^{2n}a_n-q^{2n+2}a_n+q^{2n-1}a_{n-1}+q^{4n+2}a_n,
\end{equation*}
\begin{equation*}
a_n=\frac{q^{2n-1}}{(1-q^{2n})(1+q^{2n+2})}a_{n-1}=\cdots=\frac{q^{n^2}(1+q^2)}{(q^4;q^4)_n(1+q^{2n+2})}a_0.
\end{equation*}
In the following, we show some details for the calculation of $a_0$.

$F(z)$ is the limit of the finite determinant {\small
\begin{align*}\nonumber
D_n(z):= \left| \begin{array}{ccccc} 1-q^2+zq & q^2 & & & \cdots \\
-1 & 1-q^2+zq^3 & q^2 & & \cdots  \\
& -1 & 1-q^2+zq^5 & q^2 & \cdots \\
\vdots & \vdots & \ddots & \ddots & \ddots \\
& & -1 & 1-q^2+zq^{2n-3} & q^2 \\
& & & -1 & 1-q^2+zq^{2n-1} \\
\end{array} \right|.
\end{align*}}
Expanding this determinant with respect to the last row, we get
\begin{equation}\label{A16-A20Dn}
D_n(z)=(1-q^2+zq^{2n-1})D_{n-1}(z)+q^2D_{n-2}(z), \qquad D_0(z)=1, \
D_1(z)=1-q^2+zq.
\end{equation}
Since $a_0=F(0)=\displaystyle\lim_{n \rightarrow \infty}D_n(0)$,
according to the recursion \eqref{A16-A20Dn}, we have
\begin{equation*}
D_n(0)=(1-q^2)D_{n-1}(0)+q^2D_{n-2}(0),\qquad D_0(0)=1, \
D_1(0)=1-q^2.
\end{equation*}
Thus, we get the following recursion
\begin{align*}
D_n(0)-D_{n-1}(0)&=-q^2(D_{n-1}(0)-D_{n-2}(0)) \\
&=\cdots\cdots \\
&=(-1)^{n-1}q^{2n-2}(D_1(0)-D_0(0))\\
&=(-1)^nq^{2n}.
\end{align*}
Then we have
\begin{equation*}
D_n(0)=D_{n-1}(0)+(-1)^{n}q^{2n}=\cdots=\frac{1+(-1)^nq^{2n+2}}{1+q^2}.
\end{equation*}
Finally, letting $n \rightarrow \infty$ in $D_n(0)$, we get
\begin{equation*} a_0=\lim_{n \rightarrow
\infty}D_n(0)=\frac{1}{1+q^2}.
\end{equation*}

Therefore, we have
\begin{equation*}
F(z)=\sum_{n=0}^{\infty}\frac{q^{n^2}}{(q^4;q^4)_n(1+q^{2n+2})}z^n,
\end{equation*}
and the left hand side of \eqref{A16-A20} can be expressed by
$F(q^{2m})+q^2F(q^{2m+2})$.

Due to \eqref{A16-A20Dn}, we have
\begin{equation}\label{A16-A20D}
D_{n-m}(q^{2m})=(1-q^2+q^{2n-1})D_{n-m-1}(q^{2m})+q^2D_{n-m-2}(q^{2m}).
\end{equation}
According to \eqref{A20P}, \eqref{A16Q}, and \eqref{A16-A20D}, we
notice that the sequences $\langle D_{n-m}(q^{2m})\rangle_n$,
$\langle P_n\rangle_n$, and $\langle Q_n\rangle_n$ satisfy the same
recursion. Set
\begin{equation}\label{A16-A20R}
D_{n-m}(q^{2m})=\lambda_m P_n+\mu_m Q_n.
\end{equation}
We can determine the parameters $\lambda_m$ and $\mu_m$ using the
initial conditions $D_{0}(q^{2m})=1$,
$D_{1}(q^{2m})=1-q^2+q^{2m+1}$, and the recursions \eqref{A20P} and
\eqref{A16Q}, which leads to the evaluations
\begin{equation*}
\lambda_m=\frac{Q_{m-1}}{P_mQ_{m-1}-P_{m-1}Q_m},
\end{equation*}
\begin{equation*}
\mu_m=\frac{P_{m-1}}{P_{m-1}Q_m-P_mQ_{m-1}}.
\end{equation*}
Indeed, we have
\begin{equation*}
P_mQ_{m-1}-P_{m-1}Q_m=(-1)^{m-1}q^{2m-1},
\end{equation*}
which can be proved by induction on $m$. Then we have simpler forms
for $\lambda_m$ and $\mu_m$ as follows:
\begin{equation}\label{A16-A20lm}
\lambda_m=(-1)^{m-1}q^{1-2m}Q_{m-1},\hskip 10mm
\mu_m=-(-1)^{m-1}q^{1-2m}P_{m-1}.
\end{equation}
Now setting $m\rightarrow m+1$ in \eqref{A16-A20R}, we get
\begin{equation*}
D_{n-m-1}(q^{2m+2})=\lambda_{m+1} P_n+\mu_{m+1} Q_n.
\end{equation*}
Thus, we have
\begin{align*}
\sum_{n=0}^{\infty}\frac{q^{n^2+2mn}}{(q^4;q^4)_n}&=F(q^{2m})+q^2F(q^{2m+2}) \nonumber\\
&=(\lambda_m+q^2\lambda_{m+1})P_{\infty}+(\mu_m+q^2\mu_{m+1})Q_{\infty}.
\end{align*}
According to \eqref{A16-A20lm}, we get
\begin{align*}
\lambda_m+q^2\lambda_{m+1}&=(-1)^{m}q^{1-2m}(Q_m-Q_{m-1}),\\
\mu_m+q^2\mu_{m+1}&=(-1)^{m-1}q^{1-2m}(P_m-P_{m-1}).
\end{align*}

Setting $A_m=(-1)^{m}q^{1-2m}(Q_m-Q_{m-1})$, due to \eqref{A16Q}, we
have
\begin{align*}
A_m&=xA_{m-1}+yA_{m-2}  \\
&=x(-1)^{m-1}q^{3-2m}(Q_{m-1}-Q_{m-2})+y(-1)^{m}q^{5-2m}(Q_{m-2}-Q_{m-3})  \\
&=[(-1)^{m-1}(1-q^{5-2m})x+(-1)^mq^{5-2m}y]Q_{m-2}+(-1)^{m-1}q^{5-2m}(x+y)Q_{m-3}.
\end{align*}
and
\begin{align*}
A_m&=(-1)^{m}q^{1-2m}(Q_m-Q_{m-1})  \\
&=(-1)^m(q^{5-2m}+q^{2m-3}-q^2)Q_{m-2}+(-1)^m(q^2-q^{5-2m})Q_{m-3}.
\end{align*}
Therefore, we get
\begin{equation*}
\left\{ \begin{array}{l} (-1)^{m-1}(1-q^{5-2m})x+(-1)^mq^{5-2m}y=(-1)^m(q^{5-2m}+q^{2m-3}-q^2), \\
 (-1)^{m-1}q^{5-2m}(x+y)=(-1)^m(q^2-q^{5-2m}).\end{array} \right.
\end{equation*}
Then
\begin{equation*}
\left\{ \begin{array}{l} x=-q^{2m-3}, \\
 y=1, \end{array} \right.
\end{equation*}
which means that
\begin{equation*}
A_m=-q^{2m-3}A_{m-1}+A_{m-2},\qquad A_0=1,\ A_1=0.
\end{equation*}
Similarly, setting $B_m=(-1)^{m-1}q^{1-2m}(P_m-P_{m-1})$, we obtain
\begin{equation*}
B_m=-q^{2m-3}B_{m-1}+B_{m-2},\qquad B_0=0,\ B_1=1.
\end{equation*}
Thus the above analysis has led to \eqref{A16-A20}. \qed

Setting $m=1$ and $m=0$ in \eqref{A16-A20}, we get the identities
\eqref{A16} and \eqref{A20}, respectively.

\begin{thm}\label{A29-A50Thm}We have
\begin{itemize}
\item[(1)]
\begin{align}\label{A29-A50-1}
\sum_{n=0}^{\infty}\frac{(-q;q^2)_nq^{n^2+2mn}}{(q;q)_{2n+1}}&=
\frac{q^{1-m}(q^2, q^{10},
q^{12};q^{12})_{\infty}}{(q;q^2)_{m-1}(q;q)_{\infty}}P_{m-1}
\nonumber \\
&-\frac{q^{1-m}(-q^2, -q^4,
q^6;q^6)_{\infty}(-q;q^2)_{\infty}}{(q;q^2)_{m-1}(q^2;q^2)_{\infty}}Q_{m-1},
\end{align}
where
\begin{align*}
P_m &=(1+q+q^{2m-1})P_{m-1}+(q^{2m-2}-q)P_{m-2}, \qquad P_0=1, \
P_1=1+q,\\
Q_m &=(1+q+q^{2m-1})Q_{m-1}+(q^{2m-2}-q)Q_{m-2}, \qquad Q_0=0, \
Q_1=1.
\end{align*}
\item[(2)]
\begin{equation}\label{A29-A50-2}
\sum_{n=0}^{\infty}\frac{(-q;q^2)_nq^{n^2+2mn}}{(q;q)_{2n}}=
\frac{(-q^2, -q^4,
q^6;q^6)_{\infty}(-q;q^2)_{\infty}}{(q;q^2)_m(q^2;q^2)_{\infty}}A_m
-\frac{(q^2, q^{10},
q^{12};q^{12})_{\infty}}{(q;q^2)_m(q;q)_{\infty}}B_m,
\end{equation}
where
\begin{align*}
A_m &=(1+q+q^{2m-2})A_{m-1}+(q^{2m-2}-q)A_{m-2},\qquad A_0=1, \
A_1=1,\\
B_m &=(1+q+q^{2m-2})B_{m-1}+(q^{2m-2}-q)B_{m-2},\qquad B_0=0, \
B_1=2q.
\end{align*}
\end{itemize}
\end{thm}
\pf The identities A.29 and A.50 in Slater's list are stated as
follows.

\textbf{Identity A.29 (Slater \cite{Slater}):}
\begin{equation}\label{A29}
\sum_{n=0}^{\infty}\frac{(-q;q^2)_nq^{n^2}}{(q;q)_{2n}}=
\frac{(-q^2,-q^4,q^6;q^6)_{\infty}(-q;q^2)_{\infty}}{(q^2;q^2)_{\infty}},
\end{equation}
\begin{equation}\label{A29P}
P_n=(1+q+q^{2n-1})P_{n-1}+(q^{2n-2}-q)P_{n-2}, \qquad
P_{-1}=-\frac{q}{1-q}, \ P_0=1, \ P_1=1+q.
\end{equation}

\textbf{Identity A.50 (Slater \cite{Slater}):}
\begin{equation}\label{A50}
\sum_{n=0}^{\infty}\frac{(-q;q^2)_nq^{n^2+2n}}{(q;q)_{2n+1}}=
\frac{(q^2,q^{10},q^{12};q^{12})_{\infty}}{(q;q)_{\infty}},
\end{equation}
\begin{equation}\label{A50P}
P_n=(1+q+q^{2n+1})P_{n-1}+(q^{2n}-q)P_{n-2}, \qquad P_{-1}=0, \
P_0=1, \ P_1=1+q+q^3.
\end{equation}

For the recursion \eqref{A50P}, letting $Q_n=P_{n-1}$, we get the
recursion
\begin{equation}\label{A50Q}
Q_n=(1+q+q^{2n-1})Q_{n-1}+(q^{2n-2}-q)Q_{n-2}, \qquad
Q_{-1}=\frac{1}{1-q},\ Q_0=0, \ Q_1=1.
\end{equation}
The polynomials $P_n$ in \eqref{A29P} and $Q_n$ in \eqref{A50Q}
satisfy the same recursion with different initial conditions, and
converge to the right hand sides of \eqref{A29} and \eqref{A50},
respectively.

Consider the following determinant:
\begin{equation*}
F(z):=\left| \begin{array}{ccccc} 1+q+zq & zq^2-q & & & \cdots \\
-1 & 1+q+zq^3 & zq^4-q & & \cdots  \\
& -1 & 1+q+zq^5 & zq^6-q & \cdots \\
& & \ddots & \ddots & \ddots\\
\end{array}
\right|.
\end{equation*}
Expanding the determinant with respect to the first column, we get
\begin{equation*}
F(z)=(1+q+zq)F(zq^2)+(zq^2-q)F(zq^4).
\end{equation*}
Setting
\begin{equation*}
F(z)=\sum_{n=0}^{\infty}a_nz^n,
\end{equation*}
we get, upon comparing coefficients,
\begin{equation*}
a_n=q^{2n}a_n+q^{2n+1}a_n+q^{2n-1}a_{n-1}+q^{4n-2}a_{n-1}-q^{4n+1}a_n,
\end{equation*}
\begin{equation*}
a_n=\frac{(1+q^{2n-1})q^{2n-1}}{(1-q^{2n})(1-q^{2n+1})}a_{n-1}=\cdots=\frac{(-q;q^2)_nq^{n^2}(1-q)}{(q;q)_{2n+1}}a_{0}.
\end{equation*}
Resorting to the same technique for $a_0$ in the proof of Theorem
\ref{A16-A20Thm}, we have $a_0=\frac{1}{1-q}$. Thus, we have
\begin{equation*}
F(z)=\sum_{n=0}^{\infty}\frac{(-q;q^2)_nq^{n^2}}{(q;q)_{2n+1}}z^n.
\end{equation*}
We observe that the left hand sides of \eqref{A29-A50-1} and
\eqref{A29-A50-2} can be expressed by $F(q^{2m})$ and
$F(q^{2m})-qF(q^{2m+2})$, respectively.

On the other hand, $F(z)$ is the limit of the finite determinant
{\small
\begin{align*}\nonumber
D_n(z):= \left| \begin{array}{ccccc} 1+q+zq & zq^2-q & & & \cdots \\
-1 & 1+q+zq^3 & zq^4-q & & \cdots  \\
& -1 & 1+q+zq^5 & zq^6-q & \cdots \\
\vdots & \vdots & \ddots & \ddots & \ddots \\
& & -1 & 1+q+zq^{2n-3} & zq^{2n-2}-q \\
& & & -1 & 1+q+zq^{2n-1} \\
\end{array} \right|.
\end{align*}}
Expanding this determinant with respect to the last row, we get
\begin{align*}
&D_n(z)=(1+q+zq^{2n-1})D_{n-1}(z)+(zq^{2n-2}-q)D_{n-2}(z),\\
&\qquad\qquad D_0(z)=1,\  D_1(z)=1+q+zq.
\end{align*}
Then we have
\begin{equation}\label{A29-A50D1}
D_{n-m}(q^{2m})=(1+q+q^{2n-1})D_{n-m-1}(q^{2m})+(q^{2n-2}-q)D_{n-m-2}(q^{2m}).
\end{equation}
According to \eqref{A29P}, \eqref{A50Q}, and \eqref{A29-A50D1}, we
notice that the sequences $\langle D_{n-m}(q^{2m})\rangle_n$,
$\langle P_n\rangle_n$, and $\langle Q_n\rangle_n$ satisfy the same
recursion. Set
\begin{equation*}
D_{n-m}(q^{2m})=\lambda_m P_n+\mu_m Q_n.
\end{equation*}
We can determine the parameters $\lambda_m$ and $\mu_m$ using the
initial conditions $D_{0}(q^{2m})=1$, $D_{1}(q^{2m})=1+q+q^{2m+1}$,
and the recursions \eqref{A29P} and \eqref{A50Q}, which leads to the
evaluations
\begin{equation*}
\lambda_m=\frac{Q_{m-1}}{P_mQ_{m-1}-P_{m-1}Q_m},
\end{equation*}
\begin{equation*}
\mu_m=\frac{P_{m-1}}{P_{m-1}Q_m-P_mQ_{m-1}}.
\end{equation*}
Notice that
\begin{equation*}
P_mQ_{m-1}-P_{m-1}Q_m=-q^{m-1}(q;q^2)_{m-1},
\end{equation*}
which can be proved by induction on $m$. Then we have simpler forms
for $\lambda_m$ and $\mu_m$ as follows:
\begin{equation}\label{A29-A50lm}
\lambda_m=-\frac{q^{1-m}}{(q;q^2)_{m-1}}Q_{m-1},\hskip 10mm
\mu_m=\frac{q^{1-m}}{(q;q^2)_{m-1}}P_{m-1}.
\end{equation}
Therefore, we obtain Equation \eqref{A29-A50-1}.

Meanwhile, we have
\begin{align*}
\sum_{n=0}^{\infty}\frac{(-q;q^2)_nq^{n^2+2mn}}{(q;q)_{2n}}&=F(q^{2m})-qF(q^{2m+2}) \nonumber\\
&=(\lambda_m-q\lambda_{m+1})P_{\infty}+(\mu_m-q\mu_{m+1})Q_{\infty}.
\end{align*}

According to \eqref{A29-A50lm}, we get
\begin{equation*}
\lambda_m-q\lambda_{m+1}=\frac{q^{1-m}}{(q;q^2)_m}[Q_m-(1-q^{2m-1})Q_{m-1}],
\end{equation*}
\begin{equation*}
\mu_m-q\mu_{m+1}=-\frac{q^{1-m}}{(q;q^2)_m}[P_m-(1-q^{2m-1})P_{m-1}].
\end{equation*}

Setting $A_m=q^{1-m}[Q_m-(1-q^{2m-1})Q_{m-1}]$ and
$B_m=q^{1-m}[P_m-(1-q^{2m-1})P_{m-1}]$, we get Equation
\eqref{A29-A50-2} as desired. \qed

The identities \eqref{A29} and \eqref{A50} are the special cases of
\eqref{A29-A50-2} and \eqref{A29-A50-1}, respectively.

\begin{thm}\label{A34-A36Thm}We have
\begin{equation}\label{A34-A36}
\sum_{n=0}^{\infty}\frac{(-q;q^2)_nq^{n^2+2mn}}{(q^2;q^2)_n}=
\frac{(-1)^mq^{m-m^2}Q_{m-1}}{(q,q^4,q^7;q^8)_{\infty}}-\frac{(-1)^mq^{m-m^2}P_{m-1}}{(q^3,q^4,q^5;q^8)_{\infty}},
\end{equation}
where
\begin{align*}
P_m &=(1+q^{2m-1})P_{m-1}+q^{2m-2}P_{m-2},\qquad P_{-1}=0,\ P_0=1, \
P_1=1+q,\\
Q_m &=(1+q^{2m-1})Q_{m-1}+q^{2m-2}Q_{m-2},\qquad Q_{-1}=1,\ Q_0=0, \
Q_1=1.
\end{align*}
\end{thm}
\pf We use the identities A.34 and A.36 in Slater's list to prove
the theorem.

\textbf{Identity A.34 (Slater \cite{Slater}): The analytic version
of the second G\"{o}llnitz-Gordon partition
identity.}\footnote{There is a typo in the recursion of Identity
A.34 given by Sills \cite{Sills}. In
\cite{Andrews-Knopfmacher-Paule-Prodinger}, Andrews et al. pointed
out this recursion by considering a special case of the
Al-Salam/Ismail polynomials \cite{Al-Salam-Ismail}.}

\begin{equation}\label{A34}
\sum_{n=0}^{\infty}\frac{(-q;q^2)_nq^{n^2+2n}}{(q^2;q^2)_n}=
\frac{1}{(q^3,q^4,q^5;q^8)_{\infty}},
\end{equation}
\begin{equation}\label{A34P}
P_n=(1+q^{2n+1})P_{n-1}+q^{2n}P_{n-2},\qquad P_{-1}=0,\ P_0=1, \
P_1=1+q^3.
\end{equation}

\textbf{Identity A.36 (Slater \cite{Slater}): The analytic version
of the first G\"{o}llnitz-Gordon partition identity.}

\begin{equation}\label{A36}
\sum_{n=0}^{\infty}\frac{(-q;q^2)_nq^{n^2}}{(q^2;q^2)_n}=
\frac{1}{(q,q^4,q^7;q^8)_{\infty}},
\end{equation}
\begin{equation}\label{A36P}
P_n=(1+q^{2n-1})P_{n-1}+q^{2n-2}P_{n-2},\qquad P_{-1}=0,\ P_0=1, \
P_1=1+q.
\end{equation}

For the recursion \eqref{A34P}, letting $Q_n=P_{n-1}$, we get the
recursion
\begin{equation}\label{A34Q}
Q_n=(1+q^{2n-1})Q_{n-1}+q^{2n-2}Q_{n-2},\qquad Q_{-1}=1,\ Q_0=0, \
Q_1=1.
\end{equation}
Therefore, $P_n$ in \eqref{A36P} and $Q_n$ in \eqref{A34Q} converge
to the right hand sides of \eqref{A36} and \eqref{A34},
respectively. In the following, they are used to prove this theorem.

Consider the following determinant:
\begin{equation*}
F(z):=\left| \begin{array}{ccccc} 1+zq & zq^2 & & & \cdots \\
-1 & 1+zq^3 & zq^4 & & \cdots  \\
& -1 & 1+zq^5 & zq^6 & \cdots \\
& & \ddots & \ddots & \ddots\\
\end{array}
\right|.
\end{equation*}
Expanding the determinant with respect to the first column, we get
\begin{equation*}
F(z)=(1+zq)F(zq^2)+zq^2F(zq^4).
\end{equation*}
Setting
\begin{equation*}
F(z)=\sum_{n=0}^{\infty}a_nz^n,
\end{equation*}
we get, upon comparing coefficients,
\begin{equation*}
a_n=q^{2n}a_n+q^{2n-1}a_{n-1}+q^{4n-2}a_{n-1},
\end{equation*}
\begin{equation*}
a_n=\frac{(1+q^{2n-1})q^{2n-1}}{1-q^{2n}}a_{n-1}=\cdots=\frac{(-q;q^2)_nq^{n^2}}{(q^2;q^2)_n}a_0.
\end{equation*}
Since $a_0=1$, iteration leads to
\begin{equation*}
F(z)=\sum_{n=0}^{\infty}\frac{(-q;q^2)_nq^{n^2}}{(q^2;q^2)_n}z^n,
\end{equation*}
and thus the left hand side of \eqref{A34-A36} can be expressed by
$F(q^{2m})$.

On the other hand, $F(z)$ is the limit of the finite determinant
{\small
\begin{align*}\nonumber
D_n(z):= \left| \begin{array}{ccccc} 1+zq & zq^2 & & & \cdots \\
-1 & 1+zq^3 & zq^4 & & \cdots  \\
& -1 & 1+zq^5 & zq^6 & \cdots \\
\vdots & \vdots & \ddots & \ddots & \ddots \\
& & -1 & 1+zq^{2n-3} & zq^{2n-2} \\
& & & -1 & 1+zq^{2n-1} \\
\end{array} \right|.
\end{align*}}
Expanding this determinant with respect to the last row, we get
\begin{equation*}
D_n(z)=(1+zq^{2n-1})D_{n-1}(z)+zq^{2n-2}D_{n-2}(z),\qquad D_0(z)=1,\
D_1(z)=1+zq.
\end{equation*}
Then we have
\begin{equation*}
D_{n-m}(q^{2m})=(1+q^{2n-1})D_{n-m-1}(q^{2m})+q^{2n-2}D_{n-m-2}(q^{2m}).
\end{equation*}
Therefore, we set
\begin{equation*}
D_{n-m}(q^{2m})=\lambda_m P_n+\mu_m Q_n.
\end{equation*}
Using the initial conditions $D_{0}(q^{2m})=1$ and
$D_{1}(q^{2m})=1+q^{2m+1}$, we get
\begin{equation*}
\lambda_m=\frac{Q_{m-1}}{P_mQ_{m-1}-P_{m-1}Q_m},
\end{equation*}
\begin{equation*}
\mu_m=\frac{P_{m-1}}{P_{m-1}Q_m-P_mQ_{m-1}}.
\end{equation*}
Indeed, we have
\begin{equation*}
P_mQ_{m-1}-P_{m-1}Q_m=(-1)^mq^{m^2-m},
\end{equation*}
which can be proved by induction on $m$. Then we have simpler forms
for $\lambda_m$ and $\mu_m$ as follows:
\begin{equation*}
\lambda_m=(-1)^mq^{m-m^2}Q_{m-1},\hskip 10mm
\mu_m=-(-1)^mq^{m-m^2}P_{m-1}.
\end{equation*}
Equation \eqref{A34-A36} is proved.\qed

The identities \eqref{A34} and \eqref{A36} are the special cases of
Equation \eqref{A34-A36}.

\begin{thm}\label{A38-A39Thm}We have
\begin{itemize}
\item[(1)]
\begin{align}\label{A38-A39-1}
\sum_{n=0}^{\infty}\frac{q^{2n^2+2mn}}{(q;q)_{2n+1}}&=
\frac{q^{1-m}(q^3,q^5,q^8;q^8)_{\infty}(q^2,q^{14};q^{16})_{\infty}}{(q;q^2)_{m-1}(q;q)_{\infty}}P_{m-1}
\nonumber \\
&-\frac{q^{1-m}(q,q^7,q^8;q^8)_{\infty}(q^6,q^{10};q^{16})_{\infty}}{(q;q^2)_{m-1}(q;q)_{\infty}}Q_{m-1},
\end{align}
where
\begin{align*}
P_m &=(1+q)P_{m-1}+(q^{2m-2}-q)P_{m-2},\qquad P_0=1, \ P_1=1,\\
Q_m &=(1+q)Q_{m-1}+(q^{2m-2}-q)Q_{m-2},\qquad Q_0=0, \ Q_1=1.
\end{align*}
\item[(2)]
\begin{equation}\label{A38-A39-2}
\sum_{n=0}^{\infty}\frac{q^{2n^2+2mn}}{(q;q)_{2n}}=
\frac{(q,q^7,q^8;q^8)_{\infty}(q^6,q^{10};q^{16})_{\infty}}{(q;q^2)_m(q;q)_{\infty}}A_m
-\frac{(q^3,q^5,q^8;q^8)_{\infty}(q^2,q^{14};q^{16})_{\infty}}{(q;q^2)_m(q;q)_{\infty}}B_m,
\end{equation}
where
\begin{align*}
A_m &=(1+q)A_{m-1}+(q^{2m-2}-q)A_{m-2},\qquad A_0=1, \ A_1=1,\\
B_m &=(1+q)B_{m-1}+(q^{2m-2}-q)B_{m-2},\qquad B_0=0, \ B_1=q.
\end{align*}
\end{itemize}
\end{thm}
\pf We use the following identities A.38 and A.39 in Slater's list
to prove the theorem.

\textbf{Identity A.38 (Slater \cite{Slater}):}
\begin{equation}\label{A38}
\sum_{n=0}^{\infty}\frac{q^{2n^2+2n}}{(q;q)_{2n+1}}=
\frac{(q^3,q^5,q^8;q^8)_{\infty}(q^2,q^{14};q^{16})_{\infty}}{(q;q)_{\infty}},
\end{equation}
\begin{equation}\label{A38P}
P_n=(1+q)P_{n-1}+(q^{2n}-q)P_{n-2},\qquad P_{-1}=0, \ P_0=1, \
P_1=1+q.
\end{equation}

\textbf{Identity A.39 (Jackson \cite{Jackson}):}

\begin{equation}\label{A39}
\sum_{n=0}^{\infty}\frac{q^{2n^2}}{(q;q)_{2n}}=
\frac{(q,q^7,q^8;q^8)_{\infty}(q^6,q^{10};q^{16})_{\infty}}{(q;q)_{\infty}},
\end{equation}
\begin{equation}\label{A39P}
P_n=(1+q)P_{n-1}+(q^{2n-2}-q)P_{n-2},\qquad
P_{-1}=-\frac{q}{1-q},\ P_0=1, \ P_1=1.
\end{equation}

For the recursion \eqref{A38P}, letting $Q_n=P_{n-1}$, we get the
recursion
\begin{equation}\label{A38Q}
Q_n=(1+q)Q_{n-1}+(q^{2n-2}-q)Q_{n-2},\qquad Q_{-1}=\frac{1}{1-q}, \
Q_0=0, \ Q_1=1.
\end{equation}

Therefore, $P_n$ in \eqref{A39P} and $Q_n$ in \eqref{A38Q} satisfy
the same recursion with different initial conditions, and converge
to the right hand sides of \eqref{A39} and \eqref{A38},
respectively.

Consider the following determinant:
\begin{equation*}
F(z):=\left| \begin{array}{ccccc} 1+q & zq^2-q & & & \cdots \\
-1 & 1+q & zq^4-q & & \cdots  \\
& -1 & 1+q & zq^6-q & \cdots \\
& & \ddots & \ddots & \ddots\\
\end{array}
\right|.
\end{equation*}
Expanding the determinant with respect to the first column, we get
\begin{equation*}
F(z)=(1+q)F(zq^2)+(zq^2-q)F(zq^4).
\end{equation*}
Setting
\begin{equation*}
F(z)=\sum_{n=0}^{\infty}a_nz^n,
\end{equation*}
we get, upon comparing coefficients,
\begin{equation*}
a_n=q^{2n}a_n+q^{2n+1}a_n+q^{4n-2}a_{n-1}-q^{4n+1}a_n,
\end{equation*}
\begin{equation*}
a_n=\frac{q^{4n-2}}{(1-q^{2n})(1-q^{2n+1})}a_{n-1}=\cdots=\frac{q^{2n^2}(1-q)}{(q;q)_{2n+1}}a_{0}.
\end{equation*}
Since $a_0=\frac{1}{1-q}$, iteration leads to
\begin{equation*}
F(z)=\sum_{n=0}^{\infty}\frac{q^{2n^2}}{(q;q)_{2n+1}}z^n,
\end{equation*}
and thus the left hand sides of \eqref{A38-A39-1} and
\eqref{A38-A39-2} can be expressed by $F(q^{2m})$ and
$F(q^{2m})-qF(q^{2m+2})$, respectively.

On the other hand, $F(z)$ is the limit of the finite determinant
{\small
\begin{align*}\nonumber
D_n(z):= \left| \begin{array}{ccccc} 1+q & zq^2-q & & & \cdots \\
-1 & 1+q & zq^4-q & & \cdots  \\
& -1 & 1+q & zq^6-q & \cdots \\
\vdots & \vdots & \ddots & \ddots & \ddots \\
& & -1 & 1+q & zq^{2n-2}-q \\
& & & -1 & 1+q \\
\end{array} \right|.
\end{align*}}
Expanding this determinant with respect to the last row, we get
\begin{equation*}
D_n(z)=(1+q)D_{n-1}(z)+(zq^{2n-2}-q)D_{n-2}(z), \qquad D_0(z)=1, \
D_1(z)=1+q.
\end{equation*} Then we have
\begin{equation*}
D_{n-m}(q^{2m})=(1+q)D_{n-m-1}(q^{2m})+(q^{2n-2}-q)D_{n-m-2}(q^{2m}).
\end{equation*}
Therefore, we set
\begin{equation*}
D_{n-m}(q^{2m})=\lambda_m P_n+\mu_m Q_n.
\end{equation*}
We can determine the parameters $\lambda_m$ and $\mu_m$ using the
initial conditions $D_{0}(q^{2m})=1$, $D_{1}(q^{2m})=1+q$, and the
recursions \eqref{A39P} and \eqref{A38Q}, which leads to the
evaluations
\begin{equation*}
\lambda_m=\frac{Q_{m-1}}{P_mQ_{m-1}-P_{m-1}Q_m},
\end{equation*}
\begin{equation*}
\mu_m=\frac{P_{m-1}}{P_{m-1}Q_m-P_mQ_{m-1}}.
\end{equation*}
We get
\begin{equation*}
P_mQ_{m-1}-P_{m-1}Q_m=-q^{m-1}(q;q^2)_{m-1},
\end{equation*}
which can be proved by induction on $m$. Then we have
\begin{equation}\label{A38-A39lm}
\lambda_m=-\frac{q^{1-m}}{(q;q^2)_{m-1}}Q_{m-1},\hskip 10mm
\mu_m=\frac{q^{1-m}}{(q;q^2)_{m-1}}P_{m-1}.
\end{equation}
Therefore, Equation \eqref{A38-A39-1} is proved.

Furthermore, we have
\begin{align*}
\sum_{n=0}^{\infty}\frac{q^{2n^2+2mn}}{(q;q)_{2n}}&=F(q^{2m})-q F(q^{2m+2}) \nonumber\\
&=(\lambda_m-q\lambda_{m+1})P_{\infty}+(\mu_m-q\mu_{m+1})Q_{\infty}.
\end{align*}

According to \eqref{A38-A39lm}, we get
\begin{align*}
\lambda_m-q\lambda_{m+1}&=\frac{q^{1-m}}{(q;q^2)_m}[Q_m-(1-q^{2m-1})Q_{m-1}],\\
\mu_m-q\mu_{m+1}&=-\frac{q^{1-m}}{(q;q^2)_m}[P_m-(1-q^{2m-1})P_{m-1}].
\end{align*}

Setting $A_m=q^{1-m}[Q_m-(1-q^{2m-1})Q_{m-1}]$ and
$B_m=q^{1-m}[P_m-(1-q^{2m-1})P_{m-1}]$, we obtain Equation
\eqref{A38-A39-2}.\qed

The identities \eqref{A38} and \eqref{A39} are the special cases
of \eqref{A38-A39-1} and \eqref{A38-A39-2}, respectively.

\begin{thm}\label{A79-A96Thm}We have
\begin{itemize}
\item[(1)]
\begin{align}\label{A79-A96-1}
\sum_{n=0}^{\infty}\frac{q^{n^2+2mn}}{(q;q)_{2n+1}}&=
\frac{q^{1-m}(q^4,q^6,q^{10};q^{10})_{\infty}(q^2,q^{18};q^{20})_{\infty}}{(q;q)_{\infty}}P_{m-1}
\nonumber \\
&-\frac{q^{1-m}(q^8,q^{12},q^{20};q^{20})_{\infty}(-q;q^2)_{\infty}}{(q^2;q^2)_{\infty}}Q_{m-1},
\end{align}
where
\begin{align*}
P_m &=(1+q+q^{2m-1})P_{m-1}-qP_{m-2},\qquad P_{-1}=1,\ P_0=1, \
P_1=1+q,\\
Q_m &=(1+q+q^{2m-1})Q_{m-1}-qQ_{m-2},\qquad Q_{-1}=-q^{-1},\ Q_0=0,
\ Q_1=1.
\end{align*}
\item[(2)]
\begin{equation}\label{A79-A96-2}
\sum_{n=0}^{\infty}\frac{q^{n^2+2mn}}{(q;q)_{2n}}=
\frac{(q^8,q^{12},q^{20};q^{20})_{\infty}(-q;q^2)_{\infty}}{(q^2;q^2)_{\infty}}A_m
-\frac{(q^4,q^6,q^{10};q^{10})_{\infty}(q^2,q^{18};q^{20})_{\infty}}{(q;q)_{\infty}}B_m,
\end{equation}
where
\begin{align*}
A_m &=(1+q+q^{2m-2})A_{m-1}-qA_{m-2},\qquad A_0=1, \ A_1=1,\\
B_m &=(1+q+q^{2m-2})B_{m-1}-qB_{m-2},\qquad B_0=0, \ B_1=q.
\end{align*}
\end{itemize}
\end{thm}
\pf The identities A.79 and A.96 are stated as follows.

\textbf{Identity A.79 (Rogers \cite{Rogers}):}
\begin{equation}\label{A79}
\sum_{n=0}^{\infty}\frac{q^{n^2}}{(q;q)_{2n}}=
\frac{(q^8,q^{12},q^{20};q^{20})_{\infty}(-q;q^2)_{\infty}}{(q^2;q^2)_{\infty}},
\end{equation}
\begin{equation}\label{A79P}
P_n=(1+q+q^{2n-1})P_{n-1}-qP_{n-2},\qquad P_{-1}=1,\ P_0=1, \
P_1=1+q.
\end{equation}

\textbf{Identity A.96 (Rogers \cite{Rogers}):}

\begin{equation}\label{A96}
\sum_{n=0}^{\infty}\frac{q^{n^2+2n}}{(q;q)_{2n+1}}=
\frac{(q^4,q^6,q^{10};q^{10})_{\infty}(q^2,q^{18};q^{20})_{\infty}}{(q;q)_{\infty}},
\end{equation}
\begin{equation}\label{A96P}
P_n=(1+q+q^{2n+1})P_{n-1}-qP_{n-2},\qquad P_{-1}=0,\ P_0=1, \
P_1=1+q+q^3.
\end{equation}

For the recursion \eqref{A96P}, letting $Q_n=P_{n-1}$, we get the
recursion
\begin{equation}\label{A96Q}
Q_n=(1+q+q^{2n-1})Q_{n-1}-qQ_{n-2},\qquad Q_{-1}=-q^{-1},\ Q_0=0, \
Q_1=1.
\end{equation}
The polynomials $P_n$ in \eqref{A79P} and $Q_n$ in \eqref{A96Q}
satisfy the same recursion with different initial conditions, and
converge to the right hand sides of \eqref{A79} and \eqref{A96},
respectively.

Consider the following determinant:
\begin{equation*}
F(z):=\left| \begin{array}{ccccc} 1+q+zq & -q & & & \cdots \\
-1 & 1+q+zq^3 & -q & & \cdots  \\
& -1 & 1+q+zq^5 & -q & \cdots \\
& & \ddots & \ddots & \ddots\\
\end{array}
\right|.
\end{equation*}
Expanding the determinant with respect to the first column, we get
\begin{equation*}
F(z)=(1+q+zq)F(zq^2)-qF(zq^4).
\end{equation*}
Setting
\begin{equation*}
F(z)=\sum_{n=0}^{\infty}a_nz^n,
\end{equation*}
we get, upon comparing coefficients,
\begin{equation*}
a_n=q^{2n}a_n+q^{2n+1}a_n+q^{2n-1}a_{n-1}-q^{4n+1}a_n,
\end{equation*}
\begin{equation*}
a_n=\frac{q^{2n-1}}{(1-q^{2n})(1-q^{2n+1})}a_{n-1}=\cdots=\frac{q^{n^2}(1-q)}{(q;q)_{2n+1}}a_{0}.
\end{equation*}
Since $a_0=\frac{1}{1-q}$, we have
\begin{equation*}
F(z)=\sum_{n=0}^{\infty}\frac{q^{n^2}}{(q;q)_{2n+1}}z^n,
\end{equation*}
and thus the left hand sides of \eqref{A79-A96-1} and
\eqref{A79-A96-2} can be expressed by $F(q^{2m})$ and
$F(q^{2m})-qF(q^{2m+2})$, respectively.

On the other hand, $F(z)$ is the limit of the finite determinant
{\small
\begin{align*}\nonumber
D_n(z):= \left| \begin{array}{ccccc} 1+q+zq & -q & & & \cdots \\
-1 & 1+q+zq^3 & -q & & \cdots  \\
& -1 & 1+q+zq^5 & -q & \cdots \\
\vdots & \vdots & \ddots & \ddots & \ddots \\
& & -1 & 1+q+zq^{2n-3} & -q \\
& & & -1 & 1+q+zq^{2n-1} \\
\end{array} \right|.
\end{align*}}
Expanding this determinant with respect to the last row, we get
\begin{equation*}
D_n(z)=(1+q+zq^{2n-1})D_{n-1}(z)-qD_{n-2}(z), \qquad D_0(z)=1,\
D_1(z)=1+q+zq.
\end{equation*}
Then we have
\begin{equation*}
D_{n-m}(q^{2m})=(1+q+q^{2n-1})D_{n-m-1}(q^{2m})-qD_{n-m-2}(q^{2m}).
\end{equation*}
Therefore, we set
\begin{equation*}
D_{n-m}(q^{2m})=\lambda_m P_n+\mu_m Q_n.
\end{equation*}
According to the initial conditions $D_{0}(q^{2m})=1$ and
$D_{1}(q^{2m})=1+q+q^{2m+1}$, we have
\begin{equation*}
\lambda_m=\frac{Q_{m-1}}{P_mQ_{m-1}-P_{m-1}Q_m},
\end{equation*}
\begin{equation*}
\mu_m=\frac{P_{m-1}}{P_{m-1}Q_m-P_mQ_{m-1}}.
\end{equation*}
Indeed, we have
\begin{equation*}
P_mQ_{m-1}-P_{m-1}Q_m=-q^{m-1},
\end{equation*}
which can be proved by induction on $m$. Then we have
\begin{equation}\label{A79-A96lm}
\lambda_m=-q^{1-m}Q_{m-1},\hskip 10mm \mu_m=q^{1-m}P_{m-1}.
\end{equation}
Therefore, we obtain Equation \eqref{A79-A96-1}.

Furthermore, we have
\begin{align*}
\sum_{n=0}^{\infty}\frac{q^{n^2+2mn}}{(q;q)_{2n}}&=F(q^{2m})-q F(q^{2m+2}) \nonumber\\
&=(\lambda_m-q\lambda_{m+1})P_{\infty}+(\mu_m-q\mu_{m+1})Q_{\infty}.
\end{align*}

According to \eqref{A79-A96lm}, we get
\begin{equation*}
\lambda_m-q\lambda_{m+1}=q^{1-m}(Q_m-Q_{m-1}),
\end{equation*}
\begin{equation*}
\mu_m-q\mu_{m+1}=-q^{1-m}(P_m-P_{m-1}).
\end{equation*}

Setting $A_m=q^{1-m}(Q_m-Q_{m-1})$ and $B_m=q^{1-m}(P_m-P_{m-1})$,
we obtain Equation \eqref{A79-A96-2}.\qed

The identities \eqref{A79} and \eqref{A96} are the special cases of
\eqref{A79-A96-2} and \eqref{A79-A96-1}, respectively.

\begin{thm}\label{A94-A99Thm}We have
\begin{itemize}
\item[(1)]
\begin{align}\label{A94-A99-1}
\sum_{n=0}^{\infty}\frac{q^{n^2+(2m+1)n}}{(q;q)_{2n+1}}&=
\frac{q^{-m}(q^3,q^7,q^{10};q^{10})_{\infty}(q^4,q^{16};q^{20})_{\infty}}{(q;q)_{\infty}}Q_{m-1}
\nonumber \\
&-\frac{q^{-m}(q,q^9,q^{10};q^{10})_{\infty}(q^8,q^{12};q^{20})_{\infty}}{(q;q)_{\infty}}P_{m-1},
\end{align}
where
\begin{align*}
P_m &=(1+q+q^{2m})P_{m-1}-qP_{m-2},\qquad P_{-1}=0,\ P_0=1, \
P_1=1+q+q^2,\\
Q_m &=(1+q+q^{2m})Q_{m-1}-qQ_{m-2},\qquad Q_{-1}=1,\ Q_0=1, \
Q_1=1+q^2.
\end{align*}
\item[(2)]
\begin{equation}\label{A94-A99-2}
\sum_{n=0}^{\infty}\frac{q^{n^2+(2m+1)n}}{(q;q)_{2n}}=
\frac{(q,q^9,q^{10};q^{10})_{\infty}(q^8,q^{12};q^{20})_{\infty}}{(q;q)_{\infty}}A_m
-\frac{(q^3,q^7,q^{10};q^{10})_{\infty}(q^4,q^{16};q^{20})_{\infty}}{(q;q)_{\infty}}B_m,
\end{equation}
where
\begin{align*}
A_m &=(1+q+q^{2m-1})A_{m-1}-qA_{m-2},\qquad A_0=1, \ A_1=1+q,\\
B_m &=(1+q+q^{2m-1})B_{m-1}-qB_{m-2},\qquad B_0=0, \ B_1=q.
\end{align*}
\end{itemize}
\end{thm}
\pf We state the identities A.94 and A.99 in Slater's list as
follows.

\textbf{Identity A.94 (Rogers \cite{Rogers}):}
\begin{equation}\label{A94}
\sum_{n=0}^{\infty}\frac{q^{n^2+n}}{(q;q)_{2n+1}}=
\frac{(q^3,q^7,q^{10};q^{10})_{\infty}(q^4,q^{16};q^{20})_{\infty}}{(q;q)_{\infty}},
\end{equation}
\begin{equation}\label{A94P}
P_n=(1+q+q^{2n})P_{n-1}-qP_{n-2},\qquad P_{-1}=0,\ P_0=1, \
P_1=1+q+q^2.
\end{equation}

\textbf{Identity A.99 (Rogers \cite{Rogers}):}
\begin{equation}\label{A99}
\sum_{n=0}^{\infty}\frac{q^{n^2+n}}{(q;q)_{2n}}=
\frac{(q,q^9,q^{10};q^{10})_{\infty}(q^8,q^{12};q^{20})_{\infty}}{(q;q)_{\infty}},
\end{equation}
\begin{equation}\label{A99Q}
Q_n=(1+q+q^{2n})Q_{n-1}-qQ_{n-2},\qquad Q_{-1}=1,\ Q_0=1, \
Q_1=1+q^2.
\end{equation}

Consider the following determinant:
\begin{equation*}
F(z):=\left| \begin{array}{ccccc} 1+q+zq^2 & -q & & & \cdots \\
-1 & 1+q+zq^4 & -q & & \cdots  \\
& -1 & 1+q+zq^6 & -q & \cdots \\
& & \ddots & \ddots & \ddots\\
\end{array}
\right|.
\end{equation*}
Expanding the determinant with respect to the first column, we get
\begin{equation*}
F(z)=(1+q+zq^2)F(zq^2)-qF(zq^4).
\end{equation*}
Setting
\begin{equation*}
F(z)=\sum_{n=0}^{\infty}a_nz^n,
\end{equation*}
we get, upon comparing coefficients,
\begin{equation*}
a_n=q^{2n}a_n+q^{2n+1}a_n+q^{2n}a_{n-1}-q^{4n+1}a_n,
\end{equation*}
\begin{equation*}
a_n=\frac{q^{2n}}{(1-q^{2n})(1-q^{2n+1})}a_{n-1}=\cdots=\frac{q^{n^2+n}(1-q)}{(q;q)_{2n+1}}a_{0}.
\end{equation*}
Since $a_0=\frac{1}{1-q}$, we have
\begin{equation*}
F(z)=\sum_{n=0}^{\infty}\frac{q^{n^2+n}}{(q;q)_{2n+1}}z^n,
\end{equation*}
and thus the left hand sides of \eqref{A94-A99-1} and
\eqref{A94-A99-2} can be expressed by $F(q^{2m})$ and
$F(q^{2m})-qF(q^{2m+2})$, respectively.

On the other hand, $F(z)$ is the limit of the finite determinant
{\small
\begin{align*}\nonumber
D_n(z):= \left| \begin{array}{ccccc} 1+q+zq^2 & -q & & & \cdots \\
-1 & 1+q+zq^4 & -q & & \cdots  \\
& -1 & 1+q+zq^6 & -q & \cdots \\
\vdots & \vdots & \ddots & \ddots & \ddots \\
& & -1 & 1+q+zq^{2n-2} & -q \\
& & & -1 & 1+q+zq^{2n} \\
\end{array} \right|.
\end{align*}}
Expanding this determinant with respect to the last row, we get
\begin{equation*}
D_n(z)=(1+q+zq^{2n})D_{n-1}(z)-qD_{n-2}(z),\qquad D_0(z)=1,\
D_1(z)=1+q+zq^2.
\end{equation*} Then we have
\begin{equation*}
D_{n-m}(q^{2m})=(1+q+q^{2n})D_{n-m-1}(q^{2m})-qD_{n-m-2}(q^{2m}).
\end{equation*}
Set
\begin{equation*}
D_{n-m}(q^{2m})=\lambda_m P_n+\mu_m Q_n.
\end{equation*}
Using the initial conditions $D_{0}(q^{2m})=1$ and
$D_{1}(q^{2m})=1+q+q^{2m+2}$, we get
\begin{equation*}
\lambda_m=\frac{Q_{m-1}}{P_mQ_{m-1}-P_{m-1}Q_m},
\end{equation*}
\begin{equation*}
\mu_m=\frac{P_{m-1}}{P_{m-1}Q_m-P_mQ_{m-1}}.
\end{equation*}
Indeed, we have
\begin{equation*}
P_mQ_{m-1}-P_{m-1}Q_m=q^m,
\end{equation*}
which can be proved by induction on $m$. Then we have simpler forms
for $\lambda_m$ and $\mu_m$ as follows:
\begin{equation}\label{A94-A99lm}
\lambda_m=q^{-m}Q_{m-1},\hskip 10mm \mu_m=-q^{-m}P_{m-1}.
\end{equation}
Therefore, we obtain Equation \eqref{A94-A99-1}.

Furthermore, we have
\begin{align*}
\sum_{n=0}^{\infty}\frac{q^{n^2+(2m+1)n}}{(q;q)_{2n}}&=F(q^{2m})-q F(q^{2m+2}) \nonumber\\
&=(\lambda_m-q\lambda_{m+1})P_{\infty}+(\mu_m-q\mu_{m+1})Q_{\infty}.
\end{align*}

According to \eqref{A94-A99lm}, we get
\begin{align*}
\lambda_m-q\lambda_{m+1} &=-q^{-m}(Q_m-Q_{m-1}),\\
\mu_m-q\mu_{m+1} &=q^{-m}(P_m-P_{m-1}).
\end{align*}

Setting $A_m=q^{-m}(P_m-P_{m-1})$ and $B_m=q^{-m}(Q_m-Q_{m-1})$, we
get Equation \eqref{A94-A99-2}.\qed

The identities \eqref{A94} and \eqref{A99} are the special cases of
\eqref{A94-A99-1} and \eqref{A94-A99-2}, respectively.

\begin{thm}\label{A25-A25'Thm}We have
\begin{align}\label{A25-A25'}
\sum_{n=0}^{\infty}\frac{(-q;q^2)_n q^{n^2+2mn}}{(q^4;q^4)_n}&=
\frac{(q^6;q^6)_{\infty}}{(-q^2;q^2)_{m-1}(q^4;q^4)_{\infty}(q^3,q^9;q^{12})_{\infty}}A_m
\nonumber \\
&-\frac{(q^3,q^3,q^6;q^6)_{\infty}(-q;q^2)_{\infty}}{(-q^2;q^2)_{m-1}(q^2;q^2)_{\infty}}B_m,
\end{align}
where
\begin{align*}
A_m &=-q^{2m-3}A_{m-1}+(1+q^{2m-4})A_{m-2},\qquad A_0=0, \ A_1=1,\\
B_m &=-q^{2m-3}B_{m-1}+(1+q^{2m-4})B_{m-2},\qquad B_0=-\frac{1}{2},
\ B_1=0.
\end{align*}
\end{thm}
\pf The identity A.25 in Slater's list is stated as follows:

\textbf{Identity A.25 (Slater \cite{Slater}):}
\begin{equation}\label{A25}
\sum_{n=0}^{\infty}\frac{(-q;q^2)_n q^{n^2}}{(q^4;q^4)_n}=
\frac{(q^3,q^3,q^6;q^6)_{\infty}(-q;q^2)_{\infty}}{(q^2;q^2)_{\infty}}.
\end{equation}
Sills \cite{Sills} gave the following recursion for \eqref{A25}.
\begin{equation}\label{A25P}
P_n=(1-q^2+q^{2n-1})P_{n-1}+(q^2+q^{2n-2})P_{n-2},\qquad
P_{-1}=\frac{q^2}{1+q^2},\ P_0=1, \ P_1=1+q.
\end{equation}
Recently, McLaughlin et al. \cite{McLaughlin-Sills-Zimmer} found a
partner to Equation \eqref{A25}.

\textbf{An identity (McLaughlin et al. \cite[Eq.
(2.7)]{McLaughlin-Sills-Zimmer}):}
\begin{equation}\label{A25'}
\sum_{n=0}^{\infty}\frac{(-q;q^2)_n q^{n^2+2n}}{(q^4;q^4)_n}=
\frac{(q^6;q^6)_{\infty}}{(q^4;q^4)_{\infty}(q^3,q^9;q^{12})_{\infty}}.
\end{equation}
For this identity, we also have
\begin{equation}\label{A25'Q}
Q_n=(1-q^2+q^{2n-1})Q_{n-1}+(q^2+q^{2n-2})Q_{n-2},
\end{equation}
where $P_n$ and $Q_n$ converge to the right hand sides of
\eqref{A25} and \eqref{A25'}, respecitively. The initial conditions
for $Q_n$ is given in the following analysis.

Consider the following determinant:
\begin{equation*}
F(z):=\left| \begin{array}{ccccc} 1-q^2+zq & q^2+zq^2 & & & \cdots \\
-1 & 1-q^2+zq^3 & q^2+zq^4 & & \cdots  \\
& -1 & 1-q^2+zq^5 & q^2+zq^6 & \cdots \\
& & \ddots & \ddots & \ddots\\
\end{array}
\right|.
\end{equation*}
Expanding the determinant with respect to the first column, we get
\begin{equation*}
F(z)=(1-q^2+zq)F(zq^2)+(q^2+zq^2)F(zq^4).
\end{equation*}
Setting
\begin{equation*}
F(z)=\sum_{n=0}^{\infty}a_nz^n,
\end{equation*}
we get, upon comparing coefficients,
\begin{equation*}
a_n=q^{2n}a_n-q^{2n+2}a_n+q^{2n-1}a_{n-1}+q^{4n+2}a_n+q^{4n-2}a_{n-1},
\end{equation*}
\begin{equation*}
a_n=\frac{(1+q^{2n-1})q^{2n-1}}{(1-q^{2n})(1+q^{2n+2})}a_{n-1}=\cdots=\frac{(-q;q^2)_nq^{n^2}(1+q^2)}{(q^4;q^4)_n(1+q^{2n+2})}a_0.
\end{equation*}
Since $a_0=\frac{1}{1+q^2}$, iteration leads to
\begin{equation*}
F(z)=\sum_{n=0}^{\infty}\frac{(-q;q^2)_nq^{n^2}}{(q^4;q^4)_n(1+q^{2n+2})}z^n,
\end{equation*}
and thus the left hand side of \eqref{A25-A25'} can be expressed by
$F(q^{2m})+q^2F(q^{2m+2})$.

On the other hand, $F(z)$ is the limit of the finite determinant
{\small \begin{align*}
D_n(z):= \left| \begin{array}{ccccc} 1-q^2+zq & q^2+zq^2& & & \cdots \\
-1 & 1-q^2+zq^3 & q^2+zq^4 & & \cdots \\
& -1 & 1-q^2+zq^5 & q^2+zq^6 & \cdots \\
\vdots & \vdots & \ddots & \ddots & \ddots \\
& & -1 & 1-q^2+zq^{2n-3} & q^2+zq^{2n-2} \\
& & & -1 & 1-q^2+zq^{2n-1} \\
\end{array} \right|.
\end{align*}}
Expanding this determinant with respect to the last row, we get
\begin{align*}
&D_n(z)=(1-q^2+zq^{2n-1})D_{n-1}(z)+(q^2+zq^{2n-2})D_{n-2}(z),\\
&\qquad\qquad D_0(z)=1,\ D_1(z)=1-q^2+zq.
\end{align*}
Then we have
\begin{equation*}
D_{n-m}(q^{2m})=(1-q^2+q^{2n-1})D_{n-m-1}(q^{2m})+(q^2+q^{2n-2})D_{n-m-2}(q^{2m}).
\end{equation*}
Noticing that $Q_{\infty}$ is $F(q^2)+q^2F(q^4)$, we have
\begin{equation*}
Q_n=D_{n-1}(q^2)+q^2D_{n-2}(q^4).
\end{equation*}
then we get the initial conditions for $Q_n$: $Q_0=1/2$ and $Q_1=1$.

Since the sequences $\langle D_{n-m}(q^{2m})\rangle_n$, $\langle
P_n\rangle_n$, and $\langle Q_n\rangle_n$ satisfy the same
recursion, we set
\begin{equation*}
D_{n-m}(q^{2m})=\lambda_m P_n+\mu_m Q_n.
\end{equation*}
According to the initial conditions $D_{0}(q^{2m})=1$ and
$D_{1}(q^{2m})=1-q^2+q^{2m+1}$, we have
\begin{equation*}
\lambda_m=\frac{Q_{m-1}}{P_mQ_{m-1}-P_{m-1}Q_m},
\end{equation*}
\begin{equation*}
\mu_m=\frac{P_{m-1}}{P_{m-1}Q_m-P_mQ_{m-1}}.
\end{equation*}
Indeed, we have
\begin{equation*}
P_mQ_{m-1}-P_{m-1}Q_m=(-1)^mq^{2m-2}(1-q)(-q^2;q^2)_{m-2},
\end{equation*}
which can be proved by induction on $m$.

Therefore, we have simpler forms for $\lambda_m$ and $\mu_m$ as
follows:
\begin{equation}\label{A25-A25'lm}
\lambda_m=\frac{(-1)^mq^{2-2m}}{(1-q)(-q^2;q^2)_{m-2}}Q_{m-1},\hskip
10mm \mu_m=-\frac{(-1)^mq^{2-2m}}{(1-q)(-q^2;q^2)_{m-2}}P_{m-1}.
\end{equation}
Moreover, we observe that
\begin{align*}
\sum_{n=0}^{\infty}\frac{(-q;q^2)_n q^{n^2+2mn}}{(q^4;q^4)_n}&=F(q^{2m})+q^2 F(q^{2m+2}) \nonumber\\
&=(\lambda_m+q^2\lambda_{m+1})P_{\infty}+(\mu_m+q^2\mu_{m+1})Q_{\infty}.
\end{align*}

According to \eqref{A25-A25'lm}, we get
\begin{equation*}
\lambda_m+q^2\lambda_{m+1}=-\frac{(-1)^mq^{2-2m}}{(1-q)(-q^2;q^2)_{m-1}}[Q_m-(1+q^{2m-2})Q_{m-1}],
\end{equation*}
\begin{equation*}
\mu_m+q^2\mu_{m+1}=\frac{(-1)^mq^{2-2m}}{(1-q)(-q^2;q^2)_{m-1}}[P_m-(1+q^{2m-2})P_{m-1}].
\end{equation*}
Setting
\begin{align*}
A_m&=(-1)^mq^{2-2m}[P_m-(1+q^{2m-2})P_{m-1}]/(1-q), \\
B_m&=(-1)^mq^{2-2m}[Q_m-(1+q^{2m-2})Q_{m-1}]/(1-q),
\end{align*}
we get Equation \eqref{A25-A25'}.\qed

The identities \eqref{A25} and \eqref{A25'} are the special cases of
\eqref{A25-A25'}, respectively.

\section{Generalizations of identities with four-term recursions}

In this section, we apply the determinant method to the
Rogers-Ramanujan type identities with the four-term recursions of
the polynomials which converge to the right hand sides of the
identities in \cite{Sills}. Moreover, we generalize some new
identities in recent papers
\cite{Bowman-McLaughin-Sills,McLaughlin-Sills-Zimmer}. During the
calculation, some properties of determinants are used to simplify
the identities.

Three identities are used to prove each theorem. For convenience,
we give the same recursions for the polynomials $P_n$, $Q_n$, and
$R_n$ by shifting the index of the recursions given by Sills in
\cite{Sills}, like the way we have done in the previous section,
where $P_n$, $Q_n$, and $R_n$ converge to the right hand sides of
the identities in Slater's list .

\begin{thm}\label{A31-A32-A33Thm}We have
\begin{itemize}
\item[(1)]
\begin{align}\label{A31-A32-A33-1}
\sum_{n=0}^{\infty}\frac{q^{2n^2+2mn}}{(q^2;q^2)_{n}(-q;q)_{2n+1}}&=
\frac{q^{-m}(q,q^6,q^7;q^7)_{\infty}}{(q^2;q^2)_{\infty}}A_m
+\frac{q^{-m}(q^2,q^5,q^7;q^7)_{\infty}}{(q^2;q^2)_{\infty}}B_m
\nonumber
\\&
 +\frac{q^{-m}(q^3,q^4,q^7;q^7)_{\infty}}{(q^2;q^2)_{\infty}}C_m,
\end{align}
where
\begin{align*}
&A_m=-(1+q^{2m-4})A_{m-1}+q^2A_{m-2}+q^2A_{m-3},\qquad A_0=-q, \ A_1=q, \ A_2=-q,\nonumber\\[5pt]
&B_m=-(1+q^{2m-4})B_{m-1}+q^2B_{m-2}+q^2B_{m-3},\qquad B_0=0, \ B_1=0, \ B_2=q,\nonumber\\[5pt]
&C_m=-(1+q^{2m-4})C_{m-1}+q^2C_{m-2}+q^2C_{m-3},\qquad C_0=1, \
C_1=0, \ C_2=0.
\end{align*}
\item[(2)]
\begin{equation}\label{A31-A32-A33-2}
\sum_{n=0}^{\infty}\frac{q^{2n^2+2mn}}{(q^2;q^2)_{n}(-q;q)_{2n}}=
\frac{(q,q^6,q^7;q^7)_{\infty}}{(q^2;q^2)_{\infty}}E_m
+\frac{(q^2,q^5,q^7;q^7)_{\infty}}{(q^2;q^2)_{\infty}}F_m
+\frac{(q^3,q^4,q^7;q^7)_{\infty}}{(q^2;q^2)_{\infty}}G_m,
\end{equation}
where
\begin{align*}
&E_m=-(q+q^{2m-3})E_{m-1}+E_{m-2}+qE_{m-3},\qquad E_0=0, \ E_1=0, \ E_2=q,\nonumber\\[5pt]
&F_m=-(q+q^{2m-3})F_{m-1}+F_{m-2}+qF_{m-3},\qquad F_0=0, \ F_1=1, \ F_2=-q,\nonumber\\[5pt]
&G_m=-(q+q^{2m-3})G_{m-1}+G_{m-2}+qG_{m-3},\qquad G_0=1, \ G_1=0, \
G_2=1.
\end{align*}
\end{itemize}
\end{thm}
\pf The identities A.31, A.32, and A.33 in Slater's list are stated
as follows.

\textbf{Identity A.31 (Rogers \cite{Rogers1} and Selberg
\cite{Selberg}): The third Rogers-Selberg identity.}
\begin{equation}\label{A31}
\sum_{n=0}^{\infty}\frac{q^{2n^2+2n}}{(q^2;q^2)_n(-q;q)_{2n+1}}=
\frac{(q,q^6,q^7;q^7)_{\infty}}{(q^2;q^2)_{\infty}},
\end{equation}
\begin{align}\label{A31P}
&P_n=(1-q-q^2)P_{n-1}+(q^{2n}-q^3+q^2+q)P_{n-2}+q^3P_{n-3},\nonumber
\\
&\qquad \qquad P_0=1, \ P_1=1-q, \ P_2=1-q+q^2+q^4.
\end{align}

\textbf{Identity A.32 (Rogers \cite{Rogers} and Selberg
\cite{Selberg}): The second Rogers-Selberg identity.}
\begin{equation}\label{A32}
\sum_{n=0}^{\infty}\frac{q^{2n^2+2n}}{(q^2;q^2)_n(-q;q)_{2n}}=
\frac{(q^2,q^5,q^7;q^7)_{\infty}}{(q^2;q^2)_{\infty}},
\end{equation}
\begin{align}\label{A32Q}
&Q_n=(1-q-q^2)Q_{n-1}+(q^{2n}-q^3+q^2+q)Q_{n-2}+q^3Q_{n-3},\nonumber
\\
&\qquad \qquad \qquad Q_0=1, \ Q_1=1, \ Q_2=1+q^4.
\end{align}

\textbf{Identity A.33 (Rogers \cite{Rogers} and Selberg
\cite{Selberg}): The first Rogers-Selberg identity.}
\begin{equation}\label{A33}
\sum_{n=0}^{\infty}\frac{q^{2n^2}}{(q^2;q^2)_n(-q;q)_{2n}}=
\frac{(q^3,q^4,q^7;q^7)_{\infty}}{(q^2;q^2)_{\infty}},
\end{equation}
\begin{align}\label{A33R}
&R_n=(1-q-q^2)R_{n-1}+(q^{2n}-q^3+q^2+q)R_{n-2}+q^3R_{n-3},\nonumber
\\
&\qquad \qquad R_0=1, \ R_1=1+q^2, \ R_2=1+q^2-q^3.
\end{align}
The polynomials $P_n$, $Q_n$, and $R_n$ converge to the right hand
sides of \eqref{A31}, \eqref{A32}, and \eqref{A33}, respectively.

Consider the following determinant: {\small
\begin{equation*}
F(z):=\left| \begin{array}{cccccc} 1-q-q^2 & zq^2-q^3+q^2+q & q^3 & & & \cdots \\
-1 & 1-q-q^2 & zq^4-q^3+q^2+q & q^3 & & \cdots  \\
& -1 & 1-q-q^2 & zq^6-q^3+q^2+q & q^3 & \cdots \\
& & \ddots & \ddots & \ddots & \ddots\\
\end{array}
\right|.
\end{equation*}}
Expanding the determinant with respect to the first column, we get
\begin{equation*}
F(z)=(1-q-q^2)F(zq^2)+(zq^2-q^3+q^2+q)F(zq^4)+q^3F(zq^6).
\end{equation*}
Setting
\begin{equation*}
F(z)=\sum_{n=0}^{\infty}a_nz^n,
\end{equation*}
we get, upon comparing coefficients,
\begin{equation*}
a_n=q^{2n}a_n-q^{2n+1}a_n-q^{2n+2}a_{n}+q^{4n-2}a_{n-1}-q^{4n+3}a_n+q^{4n+2}a_n+q^{4n+1}a_n
+q^{6n+3}a_n,\end{equation*}
\begin{equation*}
a_n=\frac{q^{4n-2}}{(1-q^{2n})(1+q^{2n+1})(1+q^{2n+2})}a_{n-1}=\cdots=\frac{q^{2n^2}(1+q)(1+q^2)}{(q^2;q^2)_n(-q;q)_{2n+2}}a_{0}.
\end{equation*}
Since $a_0=\frac{1}{(1+q)(1+q^2)}$, we have
\begin{equation*}
F(z)=\sum_{n=0}^{\infty}\frac{q^{2n^2}}{(q^2;q^2)_n(-q;q)_{2n+2}}z^n.
\end{equation*}
Thus we get
\begin{equation}\label{A31-A32-A33F1}
\sum_{n=0}^{\infty}\frac{q^{2n^2+2mn}}{(q^2;q^2)_{n}(-q;q)_{2n+1}}=F(q^{2m})+q^2F(q^{2m+2}),
\end{equation}
\begin{equation}\label{A31-A32-A33F2}
\sum_{n=0}^{\infty}\frac{q^{2n^2+2mn}}{(q^2;q^2)_{n}(-q;q)_{2n}}=F(q^{2m})+(q+q^2)F(q^{2m+2})+q^3F(q^{2m+4}).
\end{equation}

On the other hand, $F(z)$ is the limit of the finite determinant
{\footnotesize
\begin{align*}\nonumber
D_n(z):= \left| \begin{array}{ccccc} 1-q-q^2 & zq^2-q^3+q^2+q & q^3 & & \cdots \\
-1 & 1-q-q^2 & zq^4-q^3+q^2+q & q^3 & \cdots  \\
\vdots & \ddots & \ddots & \ddots & \ddots \\
& & -1 & 1-q-q^2 & zq^{2n-2}-q^3+q^2+q \\
& & & -1 & 1-q-q^2 \\
\end{array} \right|.
\end{align*}}
Expanding this determinant with respect to the last row, we get
\begin{align*}
&D_n(z)=(1-q-q^2)D_{n-1}(z)+(zq^{2n-2}-q^3+q^2+q)D_{n-2}(z)+q^3D_{n-3}(z),\\
&\qquad D_0(z)=1,\ D_1(z)=1-q-q^2,\ D_2(z)=1-q+q^3+q^4+zq^2.
\end{align*}
Then we have
\begin{equation*}
D_{n-m+1}(q^{2m})=(1-q-q^2)D_{n-m}(q^{2m})+(q^{2n}-q^3+q^2+q)D_{n-m-1}(q^{2m})+q^3D_{n-m-2}(q^{2m}).
\end{equation*}
Since $\langle D_{n-m+1}(q^{2m})\rangle_n$, $\langle P_n\rangle_n$,
$\langle Q_n\rangle_n$, and $\langle R_n\rangle_n$ satisfy the same
recursion, we set
\begin{equation*}
D_{n-m+1}(q^{2m})=\lambda_m P_n+\mu_m Q_n+\nu_m R_n.
\end{equation*}
Using the initial conditions $D_{0}(q^{2m})=1$,
$D_{1}(q^{2m})=1-q-q^2$, and $D_{2}(q^{2m})=1-q+q^3+q^4+q^{2m+2}$,
we have
\begin{align*}
\lambda_m&=\frac{\left| \begin{array}{ccc}
1 & Q_{m-1} & R_{m-1} \\
1-q-q^2 & Q_m & R_m \\
1-q+q^3+q^4+q^{2m+2} & Q_{m+1} & R_{m+1} \\
\end{array}\right|}{\left|\begin{array}{ccc}
P_{m-1} & Q_{m-1} & R_{m-1} \\
P_m & Q_m & R_m \\
P_{m+1} & Q_{m+1} & R_{m+1} \\
\end{array}\right|},\\
\mu_m&=\frac{\left| \begin{array}{ccc}
P_{m-1} & 1 &  R_{m-1} \\
P_m & 1-q-q^2 &  R_m \\
P_{m+1} & 1-q+q^3+q^4+q^{2m+2} & R_{m+1} \\
\end{array}\right|}{\left|\begin{array}{ccc}
P_{m-1} & Q_{m-1} & R_{m-1} \\
P_m & Q_m & R_m \\
P_{m+1} & Q_{m+1} & R_{m+1} \\
\end{array}\right|},\\
\nu_m&=\frac{\left| \begin{array}{ccc}
P_{m-1} & Q_{m-1} & 1  \\
P_m & Q_m & 1-q-q^2    \\
P_{m+1} & Q_{m+1} & 1-q+q^3+q^4+q^{2m+2} \\
\end{array}\right|}{\left|\begin{array}{ccc}
P_{m-1} & Q_{m-1} & R_{m-1} \\
P_m & Q_m & R_m \\
P_{m+1} & Q_{m+1} & R_{m+1} \\
\end{array}\right|}.
\end{align*}
Indeed, we have
\begin{equation}\label{A31-A32-A33pqr}
\left|\begin{array}{ccc}
P_{m-1} & Q_{m-1} & R_{m-1} \\
P_m & Q_m & R_m \\
P_{m+1} & Q_{m+1} & R_{m+1} \\
\end{array}\right|=-q^{3m+2}.
\end{equation}
The proof of \eqref{A31-A32-A33pqr} is by induction on $m$. The case
$m=0$ is trivial.
\begin{equation*}
\left|\begin{array}{ccc}
P_0 & Q_0 & R_0 \\
P_1 & Q_1 & R_1 \\
P_2 & Q_2 & R_2 \\
\end{array}\right|=-q^5.
\end{equation*}
The recursions \eqref{A31P}, \eqref{A32Q}, \eqref{A33R}, and some
properties of determinants are used in the following induction
step. {\small
\begin{align*}
&\quad\left|\begin{array}{ccc}
P_m & Q_m & R_m \\
P_{m+1} & Q_{m+1} & R_{m+1} \\
P_{m+2} & Q_{m+2} & R_{m+2} \\
\end{array}\right|\\
&=\left|\begin{array}{ccc}
P_m & Q_m & R_m \\
P_{m+1} & Q_{m+1} & R_{m+1} \\
(1-q-q^2)P_{m+1} & (1-q-q^2)Q_{m+1} & (1-q-q^2)R_{m+1} \\
\end{array}\right|\\
&\quad +\left|\begin{array}{ccc}
P_m & Q_m & R_m \\
P_{m+1} & Q_{m+1} & R_{m+1} \\
(q^{2m+4}-q^3+q^2+q)P_m & (q^{2m+4}-q^3+q^2+q)Q_m & (q^{2m+4}-q^3+q^2+q)R_m \\
\end{array}\right|\\
&\quad+\left|\begin{array}{ccc}
P_m & Q_m & R_m \\
P_{m+1} & Q_{m+1} & R_{m+1} \\
q^3P_{m-1} & q^3Q_{m-1} & q^3R_{m-1} \\
\end{array}\right|\\
&=q^3\left|\begin{array}{ccc}
P_{m-1} & Q_{m-1} & R_{m-1} \\
P_m & Q_m & R_m \\
P_{m+1} & Q_{m+1} & R_{m+1} \\
\end{array}\right|.
\end{align*}}
Therefore, we have simpler forms for $\lambda_m$, $\mu_m$, and
$\nu_m$ as follows:
\begin{align}\label{A31-A32-A33lambda}
\lambda_m&=-\frac{\left| \begin{array}{ccc}
1 & Q_{m-1} & R_{m-1} \\
1-q-q^2 & Q_m & R_m \\
1-q+q^3+q^4+q^{2m+2} & Q_{m+1} & R_{m+1} \\
\end{array}\right|}{q^{3m+2}},\\
\mu_m&=-\frac{\left| \begin{array}{ccc}
P_{m-1} & 1 &  R_{m-1} \\
P_m & 1-q-q^2 &  R_m \\
P_{m+1} & 1-q+q^3+q^4+q^{2m+2} & R_{m+1} \\
\end{array}\right|}{q^{3m+2}}, \nonumber\\
\nu_m&=-\frac{\left| \begin{array}{ccc}
P_{m-1} & Q_{m-1} & 1  \\
P_m & Q_m & 1-q-q^2    \\
P_{m+1} & Q_{m+1} & 1-q+q^3+q^4+q^{2m+2} \\
\end{array}\right|}{q^{3m+2}}. \nonumber
\end{align}
According to \eqref{A31-A32-A33F1} and \eqref{A31-A32-A33F2}, by
setting
\begin{equation*}
\left\{\begin{array}{l} A_m=q^m(\lambda_m+q^2\lambda_{m+1}),\\
B_m=q^m(\mu_m+q^2\mu_{m+1}),\\
C_m=q^m(\nu_m+q^2\nu_{m+1}),
\end{array}\right. \qquad \text{and} \qquad
\left\{\begin{array}{l} E_m=\lambda_m+(q+q^2)\lambda_{m+1}+q^3\lambda_{m+2},\\
F_m=\mu_m+(q+q^2)\mu_{m+1}+q^3\mu_{m+2},\\
G_m=\nu_m+(q+q^2)\nu_{m+1}+q^3\nu_{m+2},
\end{array}\right.
\end{equation*}
we have
\begin{align*}
\sum_{n=0}^{\infty}\frac{q^{2n^2+2mn}}{(q^2;q^2)_{n}(-q;q)_{2n+1}}&=
q^{-m}A_mP_{\infty}+q^{-m}B_mQ_{\infty}+q^{-m}C_mR_{\infty},\\
\sum_{n=0}^{\infty}\frac{q^{2n^2+2mn}}{(q^2;q^2)_{n}(-q;q)_{2n}}&=
E_mP_{\infty}+F_mQ_{\infty}+G_mR_{\infty}.
\end{align*}
In the following, we only present the calculation for
$A_m=q^m(\lambda_m+q^2\lambda_{m+1})$. Others are similar.

According to \eqref{A31-A32-A33lambda}, using the same technique in
the proof of \eqref{A31-A32-A33pqr}, we have
\begin{align*}
A_m&=q^m(\lambda_m+q^2\lambda_{m+1})\\
&=-\frac{\left| \begin{array}{ccc}
1 & Q_{m-1} & R_{m-1} \\
1-q & Q_m & R_m \\
1-q+q^2+q^{2m+2} & Q_{m+1} & R_{m+1} \\
\end{array}\right|}{q^{2m+2}}\\
&=-\frac{\left| \begin{array}{ccc}
0 & Q_{m-2} & R_{m-2} \\
1 & Q_{m-1} & R_{m-1} \\
1-q & Q_m & R_m \\
\end{array}\right|}{q^{2m-1}}.
\end{align*}
Then we calculate $A_{m-1}$, $A_{m-2}$, and $A_{m-3}$. Letting the
last two columns in the determinants of $A_{m-1}$, $A_{m-2}$, and
$A_{m-3}$ be the same as those of $A_m$, we set
$A_m=xA_{m-1}+yA_{m-2}+zA_{m-3}$. Solve the equation, we get
\begin{equation*}
A_m=-(1+q^{2m-4})A_{m-1}+q^2A_{m-2}+q^2A_{m-3}.
\end{equation*}
Using \eqref{A31-A32-A33lambda} and the initial conditions of $P_n$,
$Q_n$, and $R_n$, we have $A_0=-q$, $A_1=q$, and $A_2=-q$. Following
the same way, we calculate the recursions of $B_m$, $C_m$, $E_m$,
$F_m$, and $G_m$ in turn. Then we obtain \eqref{A31-A32-A33-1} and
\eqref{A31-A32-A33-2}.\qed

Notice that \eqref{A31} is a special case of \eqref{A31-A32-A33-1},
and \eqref{A32} and \eqref{A33} are the special cases of
\eqref{A31-A32-A33-2}.

\begin{thm}\label{A59-A60-A61Thm}We have
\begin{itemize}
\item[(1)]
\begin{align}\label{A59-A60-A61-1}
\sum_{n=0}^{\infty}\frac{q^{n^2+mn}}{(q;q^2)_{n+1}(q;q)_n}&=
\frac{(q^2,q^{12},q^{14};q^{14})_{\infty}}{(q;q)_{\infty}}\lambda_m
+\frac{(q^4,q^{10},q^{14};q^{14})_{\infty}}{(q;q)_{\infty}}\mu_m
\nonumber
\\&+\frac{(q^6,q^8,q^{14};q^{14})_{\infty}}{(q;q)_{\infty}}\nu_m,
\end{align}
where
\begin{align*}
&\lambda_m=(1+q^{m-3})\lambda_{m-1}+q^{-1}\lambda_{m-2}-q^{-1}\lambda_{m-3},
\qquad \lambda_0=q, \ \lambda_1=0, \ \lambda_2=1,\\
&\mu_m=(1+q^{m-3})\mu_{m-1}+q^{-1}\mu_{m-2}-q^{-1}\mu_{m-3},\qquad \mu_0=0, \ \mu_1=1, \ \mu_2=0,\\
&\nu_m=(1+q^{m-3})\nu_{m-1}+q^{-1}\nu_{m-2}-q^{-1}\nu_{m-3},\qquad
\nu_0=1, \ \nu_1=0, \ \nu_2=0.
\end{align*}
\item[(2)]
\begin{equation}\label{A59-A60-A61-2}
\sum_{n=0}^{\infty}\frac{q^{n^2+mn}}{(q;q^2)_n(q;q)_n}=
\frac{(q^2,q^{12},q^{14};q^{14})_{\infty}}{(q;q)_{\infty}}E_m
+\frac{(q^4,q^{10},q^{14};q^{14})_{\infty}}{(q;q)_{\infty}}F_m
+\frac{(q^6,q^8,q^{14};q^{14})_{\infty}}{(q;q)_{\infty}}G_m,
\end{equation}
where
\begin{align*}
&E_m=(1+q^{m-1})E_{m-1}+qE_{m-2}-qE_{m-3},\qquad E_0=0, \ E_1=-q, \
E_2=-q-q^2,\\
&F_m=(1+q^{m-1})F_{m-1}+qF_{m-2}-qF_{m-3},\qquad F_0=0, \ F_1=0, \ F_2=-q,\\
&G_m=(1+q^{m-1})G_{m-1}+qG_{m-2}-qG_{m-3},\qquad G_0=1, \ G_1=1, \
G_2=1+q.
\end{align*}
\end{itemize}
\end{thm}
\pf The identities A.59, A.60, and A.61 in slater's list are stated
as follows.

\textbf{Identity A.59 (Rogers \cite{Rogers1}):}
\begin{equation}\label{A59}
\sum_{n=0}^{\infty}\frac{q^{n^2+2n}}{(q;q^2)_{n+1}(q;q)_n}=
\frac{(q^2,q^{12},q^{14};q^{14})_{\infty}}{(q;q)_{\infty}},
\end{equation}
\begin{equation}\label{A59P}
P_n=P_{n-1}+(q+q^n)P_{n-2}-qP_{n-3},\qquad P_0=0, \ P_1=1, \ P_2=1.
\end{equation}

\textbf{Identity A.60 (Rogers \cite{Rogers1}):}
\begin{equation}\label{A60}
\sum_{n=0}^{\infty}\frac{q^{n^2+n}}{(q;q^2)_{n+1}(q;q)_n}=
\frac{(q^4,q^{10},q^{14};q^{14})_{\infty}}{(q;q)_{\infty}},
\end{equation}
\begin{equation}\label{A60Q}
Q_n=Q_{n-1}+(q+q^n)Q_{n-2}-qQ_{n-3},\qquad \ Q_0=1, \ Q_1=1, \
Q_2=1+q+q^2.
\end{equation}

\textbf{Identity A.61 (Rogers \cite{Rogers}):}
\begin{equation}\label{A61}
\sum_{n=0}^{\infty}\frac{q^{n^2}}{(q;q^2)_n(q;q)_n}=
\frac{(q^6,q^8,q^{14};q^{14})_{\infty}}{(q;q)_{\infty}},
\end{equation}
\begin{equation}\label{A61R}
R_n=R_{n-1}+(q+q^n)R_{n-2}-qR_{n-3}, \qquad R_0=1, \ R_1=1+q, \
R_2=1+q+q^2.
\end{equation}

The polynomials $P_n$, $Q_n$, and $R_n$ converge to the right hand
sides of \eqref{A59}, \eqref{A60}, and \eqref{A61}, respectively.

Consider the following determinant:
\begin{equation*}
F(z):=\left| \begin{array}{cccccc} 1 & q+zq & -q & & & \cdots \\
-1 & 1 & q+zq^2 & -q & & \cdots  \\
& -1 & 1 & q+zq^3 & -q & \cdots \\
& & \ddots & \ddots & \ddots & \ddots\\
\end{array}
\right|.
\end{equation*}
Expanding the determinant with respect to the first column, we get
\begin{equation*}
F(z)=F(zq)+(q+zq)F(zq^2)-qF(zq^3).
\end{equation*}
Setting
\begin{equation*}
F(z)=\sum_{n=0}^{\infty}a_nz^n,
\end{equation*}
we get, upon comparing coefficients,
\begin{equation*}
a_n=q^na_n+q^{2n+1}a_n+q^{2n-1}a_{n-1}-q^{3n+1}a_n,
\end{equation*}
\begin{equation*}
a_n=\frac{q^{2n-1}}{(1-q^{2n+1})(1-q^n)}a_{n-1}=\cdots=\frac{q^{n^2}(1-q)}{(q;q^2)_{n+1}(q;q)_n}a_{0}.
\end{equation*}
Since $a_0=\frac{1}{1-q}$, we have
\begin{equation*}
F(z)=\sum_{n=0}^{\infty}\frac{q^{n^2}}{(q;q^2)_{n+1}(q;q)_n}z^n.
\end{equation*}
Thus we get
\begin{equation}\label{A59-A60-A61F1}
\sum_{n=0}^{\infty}\frac{q^{n^2+mn}}{(q;q^2)_{n+1}(q;q)_n} =F(q^m),
\end{equation}
\begin{equation}\label{A59-A60-A61F2}
\sum_{n=0}^{\infty}\frac{q^{n^2+mn}}{(q;q^2)_n(q;q)_n}=F(q^m)-qF(q^{m+2}).
\end{equation}

On the other hand, $F(z)$ is the limit of the finite determinant
\begin{align*}\nonumber
D_n(z):= \left| \begin{array}{ccccc} 1 & q+zq & -q & & \cdots \\
-1 & 1 & q+zq^2 & -q & \cdots  \\
\vdots & \ddots & \ddots & \ddots & \ddots \\
& & -1 & 1 & q+zq^{n-1} \\
& & & -1 & 1 \\
\end{array} \right|.
\end{align*}
Expanding this determinant with respect to the last row, we get
\begin{align*}
&D_n(z)=D_{n-1}(z)+(q+zq^{n-1})D_{n-2}(z)-qD_{n-3}(z),\\
&\qquad D_0(z)=1,\ D_1(z)=1,\ D_2(z)=1+q+zq.
\end{align*}
Then we have
\begin{equation*}
D_{n-m+1}(q^{m})=D_{n-m}(q^m)+(q+q^n)D_{n-m-1}(q^m)-qD_{n-m-2}(q^m).
\end{equation*}
Since $\langle D_{n-m+1}(q^m)\rangle_n$, $\langle P_n\rangle_n$,
$\langle Q_n\rangle_n$, and $\langle R_n\rangle_n$ satisfy the same
recursion, we set
\begin{equation*}
D_{n-m+1}(q^m)=\lambda_m P_n+\mu_m Q_n+\nu_m R_n.
\end{equation*}
Using the initial conditions $D_{0}(q^m)=1$, $D_{1}(q^m)=1$, and
$D_{2}(q^m)=1+q+q^{m+1}$, we have
\begin{align*}
\lambda_m&=\frac{\left| \begin{array}{ccc}
1 & Q_{m-1} & R_{m-1} \\
1 & Q_m & R_m \\
1+q+q^{m+1} & Q_{m+1} & R_{m+1} \\
\end{array}\right|}{\left|\begin{array}{ccc}
P_{m-1} & Q_{m-1} & R_{m-1} \\
P_m & Q_m & R_m \\
P_{m+1} & Q_{m+1} & R_{m+1} \\
\end{array}\right|},\\
\mu_m&=\frac{\left| \begin{array}{ccc}
P_{m-1} & 1 &  R_{m-1} \\
P_m & 1 &  R_m \\
P_{m+1} & 1+q+q^{m+1} & R_{m+1} \\
\end{array}\right|}{\left|\begin{array}{ccc}
P_{m-1} & Q_{m-1} & R_{m-1} \\
P_m & Q_m & R_m \\
P_{m+1} & Q_{m+1} & R_{m+1} \\
\end{array}\right|},\\
\nu_m&=\frac{\left| \begin{array}{ccc}
P_{m-1} & Q_{m-1} & 1  \\
P_m & Q_m & 1   \\
P_{m+1} & Q_{m+1} & 1+q+q^{m+1} \\
\end{array}\right|}{\left|\begin{array}{ccc}
P_{m-1} & Q_{m-1} & R_{m-1} \\
P_m & Q_m & R_m \\
P_{m+1} & Q_{m+1} & R_{m+1} \\
\end{array}\right|},
\end{align*}
where
\begin{equation*}
\left|\begin{array}{ccc}
P_{m-1} & Q_{m-1} & R_{m-1} \\
P_m & Q_m & R_m \\
P_{m+1} & Q_{m+1} & R_{m+1} \\
\end{array}\right|=(-1)^{m-1}q^m,
\end{equation*}
which can be proved by induction on $m$. Therefore, we have simpler
forms for $\lambda_m$, $\mu_m$, and $\nu_m$ as follows:
\begin{align*}
\lambda_m&=\frac{(-1)^m}{q^{m-1}}\left| \begin{array}{ccc}
0 & Q_{m-2} & R_{m-2} \\
1 & Q_{m-1} & R_{m-1} \\
1 & Q_m & R_m \\
\end{array}\right|,\\
\mu_m&=\frac{(-1)^m}{q^{m-1}}\left| \begin{array}{ccc}
P_{m-2} & 0 & R_{m-2} \\
P_{m-1} & 1 & R_{m-1} \\
P_m & 1 & R_m \\
\end{array}\right|,\\
\nu_m&=\frac{(-1)^m}{q^{m-1}}\left| \begin{array}{ccc}
P_{m-2} & Q_{m-2} & 0 \\
P_{m-1} & Q_{m-1} & 1 \\
P_m & Q_m & 1\\
\end{array}\right|.
\end{align*}
According to \eqref{A59-A60-A61F1} and \eqref{A59-A60-A61F2}, by
setting
\begin{equation*}
\left\{\begin{array}{l} E_m=\lambda_m-q\lambda_{m+2},\\
F_m=\mu_m-q\mu_{m+2},\\
G_m=\nu_m-q\nu_{m+2},
\end{array}\right.
\end{equation*}
we have
\begin{align*}
\sum_{n=0}^{\infty}\frac{q^{n^2+mn}}{(q;q^2)_{n+1}(q;q)_n}&=
\lambda_mP_{\infty}+\mu_mQ_{\infty}+\nu_mR_{\infty},\\
\sum_{n=0}^{\infty}\frac{q^{n^2+mn}}{(q;q^2)_n(q;q)_n}&=
E_mP_{\infty}+F_mQ_{\infty}+G_mR_{\infty}.
\end{align*}
Letting the last two columns in the determinants of $\lambda_{m-1}$,
$\lambda_{m-2}$, $\lambda_{m-3}$ be the same as those of
$\lambda_m$, we find a linear equation
\begin{equation*}
\lambda_m=(1+q^{m-3})\lambda_{m-1}+q^{-1}\lambda_{m-2}-q^{-1}\lambda_{m-3}.
\end{equation*}
Using the initial conditions of $P_n$, $Q_n$, and $R_n$, we have
\begin{equation*}
\lambda_0=q, \ \lambda_1=0, \ \lambda_2=1.
\end{equation*}
Proceeding in the same way, we get the recursions of $\mu_m$,
$\nu_m$, $E_m$, $F_m$, and $G_m$. Therefore, we obtain
\eqref{A59-A60-A61-1} and \eqref{A59-A60-A61-2}. \qed

The identities \eqref{A59} and \eqref{A60} are the special cases of
\eqref{A59-A60-A61-1}, and the identity \eqref{A61} is a special
case of \eqref{A59-A60-A61-2}.

\begin{thm}\label{A80-A81-A82Thm}We have
\begin{itemize}
\item[(1)]
\begin{align}\label{A80-A81-A82-1}
\sum_{n=0}^{\infty}\frac{q^{n(n+2m+1)/2}}{(q;q^2)_{n+1}(q;q)_n}&=
\frac{q^{-m}(q^2,q^5,q^7;q^7)_{\infty}(q^3,q^{11};q^{14})_{\infty}(-q;q)_{\infty}}{(q;q)_{\infty}}A_m
\nonumber\\
&+\frac{q^{-m}(q,q^6,q^7;q^7)_{\infty}(q^5,q^9;q^{14})_{\infty}(-q;q)_{\infty}}{(q;q)_{\infty}}B_m
\nonumber\\
&+\frac{q^{-m}(q^3,q^4,q^7;q^7)_{\infty}(q,q^{13};q^{14})_{\infty}(-q;q)_{\infty}}{(q;q)_{\infty}}C_m,
\end{align}
where
\begin{align*}
&A_m=qA_{m-1}+(q+q^{m-1})A_{m-2}-q^2A_{m-3},\qquad A_0=1, \ A_1=0, \ A_2=q,\\
&B_m=qB_{m-1}+(q+q^{m-1})B_{m-2}-q^2B_{m-3},\qquad B_0=0, \ B_1=0, \
B_2=-q,\\
&C_m=qC_{m-1}+(q+q^{m-1})C_{m-2}-q^2C_{m-3},\qquad C_0=0, \ C_1=q, \
C_2=0.
\end{align*}
\item[(2)]
\begin{align}\label{A80-A81-A82-2}
\sum_{n=0}^{\infty}\frac{q^{n(n+2m+1)/2}}{(q;q^2)_n(q;q)_n}&=
\frac{(q^2,q^5,q^7;q^7)_{\infty}(q^3,q^{11};q^{14})_{\infty}(-q;q)_{\infty}}{(q;q)_{\infty}}E_m
\nonumber\\
&+\frac{(q,q^6,q^7;q^7)_{\infty}(q^5,q^9;q^{14})_{\infty}(-q;q)_{\infty}}{(q;q)_{\infty}}F_m
\nonumber\\
&+\frac{(q^3,q^4,q^7;q^7)_{\infty}(q,q^{13};q^{14})_{\infty}(-q;q)_{\infty}}{(q;q)_{\infty}}G_m,
\end{align}
where
\begin{align*}
&E_m=E_{m-1}+(q+q^{m-1})E_{m-2}-qE_{m-3},\qquad E_0=0, \ E_1=0, \
E_2=-q,\\
&F_m=F_{m-1}+(q+q^{m-1})F_{m-2}-qF_{m-3},\qquad F_0=1, \ F_1=1, \
F_2=1+q,\\
&G_m=G_{m-1}+(q+q^{m-1})G_{m-2}-qG_{m-3},\qquad G_0=0, \ G_1=-q, \
G_2=-q.
\end{align*}
\end{itemize}
\end{thm}
\pf The identities A.80, A.81, and A.82 are stated as follows.

\textbf{Identity A.80 (Rogers \cite{Rogers1}):}
\begin{equation}\label{A80}
\sum_{n=0}^{\infty}\frac{q^{n(n+1)/2}}{(q;q^2)_{n+1}(q;q)_n}=
\frac{(q^2,q^5,q^7;q^7)_{\infty}(q^3,q^{11};q^{14})_{\infty}(-q;q)_{\infty}}{(q;q)_{\infty}},
\end{equation}
\begin{equation}\label{A80P}
P_n=(1+q^n)P_{n-1}+qP_{n-2}-qP_{n-3}, \qquad P_0=1, \ P_1=1+q, \
P_2=1+2q+q^2+q^3.
\end{equation}

\textbf{Identity A.81 (Rogers \cite{Rogers1}):}
\begin{equation}\label{A81}
\sum_{n=0}^{\infty}\frac{q^{n(n+1)/2}}{(q;q^2)_n(q;q)_n}=
\frac{(q,q^6,q^7;q^7)_{\infty}(q^5,q^9;q^{14})_{\infty}(-q;q)_{\infty}}{(q;q)_{\infty}},
\end{equation}
\begin{equation}\label{A81Q}
Q_n=(1+q^n)Q_{n-1}+qQ_{n-2}-qQ_{n-3},\qquad Q_0=1, \ Q_1=1+q, \
Q_2=1+q+q^2+q^3.
\end{equation}

\textbf{Identity A.82 (Rogers \cite{Rogers1}):}
\begin{equation}\label{A82}
\sum_{n=0}^{\infty}\frac{q^{n(n+3)/2}}{(q;q^2)_{n+1}(q;q)_n}=
\frac{(q^3,q^4,q^7;q^7)_{\infty}(q,q^{13};q^{14})_{\infty}(-q;q)_{\infty}}{(q;q)_{\infty}},
\end{equation}
\begin{equation}\label{A82R}
R_n=(1+q^n)R_{n-1}+qR_{n-2}-qR_{n-3},\qquad R_0=0, \ R_1=1, \
R_2=1+q^2.
\end{equation}

The polynomials $P_n$, $Q_n$, and $R_n$ converge to the right hand
sides of \eqref{A80}, \eqref{A81}, and \eqref{A82}, respectively.

Consider the following determinant:
\begin{equation*}
F(z):=\left| \begin{array}{cccccc} 1+zq & q & -q & & & \cdots \\
-1 & 1+zq^2 & q & -q & & \cdots  \\
& -1 & 1+zq^3 & q & -q & \cdots \\
& & \ddots & \ddots & \ddots & \ddots\\
\end{array}
\right|.
\end{equation*}
Expanding the determinant with respect to the first column, we get
\begin{equation*}
F(z)=(1+zq)F(zq)+qF(zq^2)-qF(zq^3).
\end{equation*}
Setting
\begin{equation*}
F(z)=\sum_{n=0}^{\infty}a_nz^n,
\end{equation*}
we get, upon comparing coefficients,
\begin{equation*}
a_n=q^na_n+q^na_{n-1}+q^{2n+1}a_n-q^{3n+1}a_n,
\end{equation*}
\begin{equation*}
a_n=\frac{q^n}{(1-q^{2n+1})(1-q^n)}a_{n-1}=\cdots=\frac{q^{(n^2+n)/2}(1-q)}{(q;q^2)_{n+1}(q;q)_n}a_{0}.
\end{equation*}
Since $a_0=\frac{1}{1-q}$, we have
\begin{equation*}
F(z)=\sum_{n=0}^{\infty}\frac{q^{(n^2+n)/2}}{(q;q^2)_{n+1}(q;q)_n}z^n.
\end{equation*}
Thus we get
\begin{equation}\label{A80-A81-A82F1}
\sum_{n=0}^{\infty}\frac{q^{n(n+2m+1)/2}}{(q;q^2)_{n+1}(q;q)_n}=F(q^m),
\end{equation}
\begin{equation}\label{A80-A81-A82F2}
\sum_{n=0}^{\infty}\frac{q^{n(n+2m+1)/2}}{(q;q^2)_n(q;q)_n}=F(q^m)-qF(q^{m+2}).
\end{equation}

On the other hand, $F(z)$ is the limit of the finite determinant
\begin{align*}\nonumber
D_n(z):= \left| \begin{array}{ccccc} 1+zq & q & -q & & \cdots \\
-1 & 1+zq^2 & q & -q & \cdots  \\
\vdots & \ddots & \ddots & \ddots & \ddots \\
& & -1 & 1+zq^{n-1} & q \\
& & & -1 & 1+zq^n \\
\end{array} \right|.
\end{align*}
Expanding this determinant with respect to the last row, we get
\begin{align*}
&\quad D_n(z)=(1+zq^n)D_{n-1}(z)+qD_{n-2}(z)-qD_{n-3}(z),\\
&D_0(z)=1,\ D_1(z)=1+zq,\ D_2(z)=1+q+zq+zq^2+z^2q^3.
\end{align*}
Then we have
\begin{equation*}
D_{n-m}(q^{m})=(1+q^n)D_{n-m-1}(q^m)+qD_{n-m-2}(q^m)-qD_{n-m-3}(q^m).
\end{equation*}
Since $\langle D_{n-m}(q^m)\rangle_n$, $\langle P_n\rangle_n$,
$\langle Q_n\rangle_n$, and $\langle R_n\rangle_n$ satisfy the same
recursion, we set
\begin{equation*}
D_{n-m}(q^m)=\lambda_m P_n+\mu_m Q_n+\nu_m R_n.
\end{equation*}
Using the initial conditions $D_{0}(q^m)=1$, $D_{1}(q^m)=1+q^{m+1}$,
and $D_{2}(q^m)=1+q+q^{m+1}+q^{m+2}+q^{2m+3}$, we have
\begin{align*}
\lambda_m&=\frac{\left| \begin{array}{ccc}
0 & Q_{m-1} & R_{m-1} \\
1 & Q_m & R_m \\
1+q^{m+1} & Q_{m+1} & R_{m+1} \\
\end{array}\right|}{\left|\begin{array}{ccc}
P_{m-1} & Q_{m-1} & R_{m-1} \\
P_m & Q_m & R_m \\
P_{m+1} & Q_{m+1} & R_{m+1} \\
\end{array}\right|},\\
\mu_m&=\frac{\left| \begin{array}{ccc}
P_{m-1} & 0 &  R_{m-1} \\
P_m & 1 &  R_m \\
P_{m+1} & 1+q^{m+1} & R_{m+1} \\
\end{array}\right|}{\left|\begin{array}{ccc}
P_{m-1} & Q_{m-1} & R_{m-1} \\
P_m & Q_m & R_m \\
P_{m+1} & Q_{m+1} & R_{m+1} \\
\end{array}\right|},\\
\nu_m&=\frac{\left| \begin{array}{ccc}
P_{m-1} & Q_{m-1} & 0  \\
P_m & Q_m & 1   \\
P_{m+1} & Q_{m+1} & 1+q^{m+1} \\
\end{array}\right|}{\left|\begin{array}{ccc}
P_{m-1} & Q_{m-1} & R_{m-1} \\
P_m & Q_m & R_m \\
P_{m+1} & Q_{m+1} & R_{m+1} \\
\end{array}\right|},
\end{align*}
where
\begin{equation*}
\left|\begin{array}{ccc}
P_{m-1} & Q_{m-1} & R_{m-1} \\
P_m & Q_m & R_m \\
P_{m+1} & Q_{m+1} & R_{m+1} \\
\end{array}\right|=(-1)^{m-1}q^m,
\end{equation*}
which can be proved by induction on $m$. Therefore, we have simpler
forms for $\lambda_m$, $\mu_m$, and $\nu_m$ as follows:
\begin{align*}
\lambda_m&=\frac{(-1)^m}{q^{m-1}}\left| \begin{array}{ccc}
0 & Q_{m-2} & R_{m-2} \\
0 & Q_{m-1} & R_{m-1} \\
1 & Q_m & R_m \\
\end{array}\right|,\\
\mu_m&=\frac{(-1)^m}{q^{m-1}}\left| \begin{array}{ccc}
P_{m-2} & 0 & R_{m-2} \\
P_{m-1} & 0 & R_{m-1} \\
P_m & 1 & R_m \\
\end{array}\right|,\\
\nu_m&=\frac{(-1)^m}{q^{m-1}}\left| \begin{array}{ccc}
P_{m-2} & Q_{m-2} & 0 \\
P_{m-1} & Q_{m-1} & 0 \\
P_m & Q_m & 1\\
\end{array}\right|.
\end{align*}
According to \eqref{A80-A81-A82F1} and \eqref{A80-A81-A82F2}, by
setting
\begin{equation*}
\left\{\begin{array}{l} A_m=q^m\lambda_m,\\
B_m=q^m\mu_m,\\
C_m=q^m\nu_m,
\end{array}\right. \qquad \text{and} \qquad
\left\{\begin{array}{l} E_m=\lambda_m-q\lambda_{m+2},\\
F_m=\mu_m-q\mu_{m+2},\\
G_m=\nu_m-q\nu_{m+2},
\end{array}\right.
\end{equation*}
we have
\begin{align*}
\sum_{n=0}^{\infty}\frac{q^{n(n+2m+1)/2}}{(q;q^2)_{n+1}(q;q)_n}&=
q^{-m}A_mP_{\infty}+q^{-m}B_mQ_{\infty}+q^{-m}C_mR_{\infty},\\
\sum_{n=0}^{\infty}\frac{q^{n(n+2m+1)/2}}{(q;q^2)_n(q;q)_n}&=
E_mP_{\infty}+F_mQ_{\infty}+G_mR_{\infty}.
\end{align*}
Since
\begin{equation*}
A_m=(-1)^m\left| \begin{array}{ccc}
0 & Q_{m-2} & R_{m-2} \\
0 & Q_{m-1} & R_{m-1} \\
q & Q_m & R_m \\
\end{array}\right|,
\end{equation*}
by letting the last two columns in the determinants of $A_{m-1}$,
$A_{m-2}$, and $A_{m-3}$ be the same as those of $A_m$, we find a
linear equation
\begin{equation*}
A_m=qA_{m-1}+(q+q^{m-1})A_{m-2}-q^2A_{m-3}.
\end{equation*}
Using the initial conditions of $P_n$, $Q_n$, and $R_n$, we have
\begin{equation*}
A_0=1, \ A_1=0, \ A_2=q.
\end{equation*}
Proceeding in the same way, we get the recursions of $B_m$, $C_m$,
$E_m$, $F_m$, and $G_m$. Therefore, we obtain \eqref{A80-A81-A82-1}
and \eqref{A80-A81-A82-2}. \qed

The identities \eqref{A80} and \eqref{A82} are the special cases of
\eqref{A80-A81-A82-1}, and \eqref{A81} is a special case of
\eqref{A80-A81-A82-2}.

\begin{thm}\label{A117-A118-A119Thm}We have
\begin{itemize}
\item[(1)]
\begin{align}\label{A117-A118-A119-1}
\sum_{n=0}^{\infty}\frac{(-q;q^2)_{n+1}q^{n^2+2mn}}{(q;q^2)_{2n+1}}&=
\frac{q^{-m}(q^3,q^{11},q^{14};q^{14})_{\infty}(q^8,q^{20};q^{28})_{\infty}(-q;q^2)_{\infty}}{(q^2;q^2)_{\infty}}A_m
\nonumber\\
&+\frac{q^{-m}(q,q^{13},q^{14};q^{14})_{\infty}(q^{12},q^{16};q^{28})_{\infty}(-q;q^2)_{\infty}}{(q^2;q^2)_{\infty}}B_m
\nonumber\\
&+\frac{q^{-m}(q^5,q^9,q^{14};q^{14})_{\infty}(q^4,q^{24};q^{28})_{\infty}(-q;q^2)_{\infty}}{(q^2;q^2)_{\infty}}C_m,
\end{align}
where
\begin{align*}
&A_m=A_{m-1}+(q^2+q^{2m-4})A_{m-2}-q^2A_{m-3},\qquad A_0=1, \ A_1=0, \ A_2=0,\\
&B_m=B_{m-1}+(q^2+q^{2m-4})B_{m-2}-q^2B_{m-3},\qquad B_0=0, \ B_1=0, \ B_2=-q,\\
&C_m=C_{m-1}+(q^2+q^{2m-4})C_{m-2}-q^2C_{m-3},\qquad C_0=q, \ C_1=q,
\ C_2=q.
\end{align*}
\item[(2)]
\begin{align}\label{A117-A118-A119-2}
\sum_{n=0}^{\infty}\frac{(-q;q^2)_nq^{n^2+2mn}}{(q;q^2)_{2n}}&=
\frac{(q^3,q^{11},q^{14};q^{14})_{\infty}(q^8,q^{20};q^{28})_{\infty}(-q;q^2)_{\infty}}{(q^2;q^2)_{\infty}}E_m
\nonumber\\
&+\frac{(q,q^{13},q^{14};q^{14})_{\infty}(q^{12},q^{16};q^{28})_{\infty}(-q;q^2)_{\infty}}{(q^2;q^2)_{\infty}}F_m
\nonumber\\
&+\frac{(q^5,q^9,q^{14};q^{14})_{\infty}(q^4,q^{24};q^{28})_{\infty}(-q;q^2)_{\infty}}{(q^2;q^2)_{\infty}}G_m,
\end{align}
where
\begin{align*}
&E_m=qE_{m-1}+(1+q^{2m-4})E_{m-2}-qE_{m-3},\qquad E_0=1, \ E_1=0, \ E_2=1,\\
&F_m=qF_{m-1}+(1+q^{2m-4})F_{m-2}-qF_{m-3},\qquad F_0=0, \ F_1=1, \ F_2=0,\\
&G_m=qG_{m-1}+(1+q^{2m-4})G_{m-2}-qG_{m-3},\qquad G_0=0, \ G_1=0, \
G_2=-q.
\end{align*}
\end{itemize}
\end{thm}
\pf We give the identities A.117, A.118, and A.119 as follows.

\textbf{Identity A.117 (Slater \cite{Slater}):}
\begin{equation}\label{A117}
\sum_{n=0}^{\infty}\frac{(-q;q^2)_nq^{n^2}}{(q^2;q^2)_{2n}}=
\frac{(q^3,q^{11},q^{14};q^{14})_{\infty}(q^8,q^{20};q^{28})_{\infty}(-q;q^2)_{\infty}}{(q^2;q^2)_{\infty}},
\end{equation}
\begin{align}\label{A117P}
&P_n=(1+q-q^2+q^{2n-1})P_{n-1}+(q^3+q^2-q)P_{n-2}-q^3P_{n-3},\nonumber
\\
&\qquad P_0=1, \ P_1=1+q, \ P_2=1+q+q^2+q^4.
\end{align}

\textbf{Identity A.118 (Slater \cite{Slater}):}
\begin{equation}\label{A118}
\sum_{n=0}^{\infty}\frac{(-q;q^2)_nq^{n^2+2n}}{(q^2;q^2)_{2n}}=
\frac{(q,q^{13},q^{14};q^{14})_{\infty}(q^{12},q^{16};q^{28})_{\infty}(-q;q^2)_{\infty}}{(q^2;q^2)_{\infty}},
\end{equation}
\begin{align}\label{A118Q}
&Q_n=(1+q-q^2+q^{2n-1})Q_{n-1}+(q^3+q^2-q)Q_{n-2}-q^3Q_{n-3},\nonumber
\\
&\qquad \qquad Q_0=1, \ Q_1=1, \ Q_2=1+q^3.
\end{align}

\textbf{Identity A.119 (Slater \cite{Slater}):}
\begin{equation}\label{A119}
\sum_{n=0}^{\infty}\frac{(-q;q^2)_{n+1}q^{n^2+2n}}{(q^2;q^2)_{2n+1}}=
\frac{(q^5,q^9,q^{14};q^{14})_{\infty}(q^4,q^{24};q^{28})_{\infty}(-q;q^2)_{\infty}}{(q^2;q^2)_{\infty}},
\end{equation}
\begin{align}\label{A119R}
&R_n=(1+q-q^2+q^{2n-1})R_{n-1}+(q^3+q^2-q)R_{n-2}-q^3R_{n-3},\nonumber
\\
&\qquad \qquad R_0=0, \ R_1=1, \ R_2=1+q+q^3.
\end{align}
The polynomials $P_n$, $Q_n$, and $R_n$ converge to the right hand
sides of \eqref{A117}, \eqref{A118}, and \eqref{A119}, respectively.

Consider the following determinant:
\begin{equation*}
F(z):=\left| \begin{array}{cccccc} 1+q-q^2+zq & q^3+q^2-q & -q^3 & & & \cdots \\
-1 & 1+q-q^2+zq^3 & q^3+q^2-q & -q^3 & & \cdots  \\
& -1 & 1+q-q^2+zq^5 & q^3+q^2-q & -q^3 & \cdots \\
& & \ddots & \ddots & \ddots & \ddots\\
\end{array}
\right|.
\end{equation*}
Expanding the determinant with respect to the first column, we get
\begin{equation*}
F(z)=(1+q-q^2+zq)F(zq^2)+(q^3+q^2-q)F(zq^4)-q^3F(zq^6).
\end{equation*}
Setting
\begin{equation*}
F(z)=\sum_{n=0}^{\infty}a_nz^n,
\end{equation*}
we get, upon comparing coefficients,
\begin{equation*}
a_n=\frac{q^{2n-1}}{(1-q^{2n})(1-q^{2n+1})(1+q^{2n+2})}a_{n-1}=\cdots=
\frac{q^{n^2}(1-q)(1+q^2)}{(q^2;q^2)_n(q;q^2)_{n+1}(-q^2;q^2)_{n+1}}a_{0}.
\end{equation*}
Since $a_0=\frac{1}{(1-q)(1+q^2)}$, using some calculations of the
$q$-shifted factorial, we have
\begin{equation*}
F(z)=\sum_{n=0}^{\infty}\frac{(-q;q^2)_{n+1}q^{n^2}}{(q^2;q^2)_{2n+1}(1+q^{2n+2})}z^n.
\end{equation*}
Thus we get
\begin{equation}\label{A117-A118-A119F1}
\sum_{n=0}^{\infty}\frac{(-q;q^2)_{n+1}q^{n^2+2mn}}{(q;q^2)_{2n+1}}=F(q^{2m})+q^2F(q^{2m+2}),
\end{equation}
\begin{equation}\label{A117-A118-A119F2}
\sum_{n=0}^{\infty}\frac{(-q;q^2)_nq^{n^2+2mn}}{(q;q^2)_{2n}}=F(q^{2m})+(q^2-q)F(q^{2m+2})-q^3F(q^{2m+4}).
\end{equation}

On the other hand, $F(z)$ is the limit of the finite determinant
{\scriptsize
\begin{align*}\nonumber
D_n(z):= \left| \begin{array}{ccccc} 1+q-q^2+zq & q^3+q^2-q & -q^3 & & \cdots \\
-1 & 1+q-q^2+zq^3 & q^3+q^2-q & -q^3 & \cdots  \\
\vdots & \ddots & \ddots & \ddots & \ddots \\
& & -1 & 1+q-q^2+zq^{2n-3} & q^3+q^2-q \\
& & & -1 & 1+q-q^2+zq^{2n-1} \\
\end{array} \right|.
\end{align*}}
Expanding this determinant with respect to the last row, we get
\begin{align*}
&D_n(z)=(1+q-q^2+zq^{2n-1})D_{n-1}(z)+(q^3+q^2-q)D_{n-2}(z)-q^3D_{n-3}(z),\\
&\qquad\qquad D_0(z)=1,\ D_1(z)=1+q-q^2+zq,\\
& \qquad D_2(z)=1+q-q^3+q^4+zq+zq^2+zq^4-zq^5+z^2q^4.
\end{align*}
Then we have{\small
\begin{equation*}
D_{n-m}(q^{2m})=(1+q-q^2+q^{2n-1})D_{n-m-1}(q^{2m})+(q^3+q^2-q)D_{n-m-2}(q^{2m})-q^3D_{n-m-3}(q^{2m}).
\end{equation*}}
Since $\langle D_{n-m}(q^{2m})\rangle_n$, $\langle P_n\rangle_n$,
$\langle Q_n\rangle_n$, and $\langle R_n\rangle_n$ satisfy the same
recursion, we set
\begin{equation*}
D_{n-m}(q^{2m})=\lambda_m P_n+\mu_m Q_n+\nu_m R_n.
\end{equation*}
Using the initial conditions $D_{-1}(q^{2m})=0$, $D_{0}(q^{2m})=1$,
and $D_{1}(q^{2m})=1+q-q^2+q^{2m+1}$, we have
\begin{align*}
\lambda_m&=\frac{\left| \begin{array}{ccc}
0 & Q_{m-1} & R_{m-1} \\
1 & Q_m & R_m \\
1+q-q^2+q^{2m+1} & Q_{m+1} & R_{m+1} \\
\end{array}\right|}{\left|\begin{array}{ccc}
P_{m-1} & Q_{m-1} & R_{m-1} \\
P_m & Q_m & R_m \\
P_{m+1} & Q_{m+1} & R_{m+1} \\
\end{array}\right|},\\
\mu_m&=\frac{\left| \begin{array}{ccc}
P_{m-1} & 0 &  R_{m-1} \\
P_m & 1 &  R_m \\
P_{m+1} & 1+q-q^2+q^{2m+1} & R_{m+1} \\
\end{array}\right|}{\left|\begin{array}{ccc}
P_{m-1} & Q_{m-1} & R_{m-1} \\
P_m & Q_m & R_m \\
P_{m+1} & Q_{m+1} & R_{m+1} \\
\end{array}\right|},\\
\nu_m&=\frac{\left| \begin{array}{ccc}
P_{m-1} & Q_{m-1} & 0  \\
P_m & Q_m & 1   \\
P_{m+1} & Q_{m+1} & 1+q-q^2+q^{2m+1} \\
\end{array}\right|}{\left|\begin{array}{ccc}
P_{m-1} & Q_{m-1} & R_{m-1} \\
P_m & Q_m & R_m \\
P_{m+1} & Q_{m+1} & R_{m+1} \\
\end{array}\right|},
\end{align*}
where
\begin{equation*}
\left|\begin{array}{ccc}
P_{m-1} & Q_{m-1} & R_{m-1} \\
P_m & Q_m & R_m \\
P_{m+1} & Q_{m+1} & R_{m+1} \\
\end{array}\right|=(-1)^mq^{3m},
\end{equation*}
which can be proved by induction on $m$. Therefore, we have simpler
forms for $\lambda_m$, $\mu_m$, and $\nu_m$ as follows:
\begin{align*}
\lambda_m&=\frac{(-1)^{m-1}}{q^{3m-3}}\left| \begin{array}{ccc}
0 & Q_{m-2} & R_{m-2} \\
0 & Q_{m-1} & R_{m-1} \\
1 & Q_m & R_m \\
\end{array}\right|,\\
\mu_m&=\frac{(-1)^{m-1}}{q^{3m-3}}\left| \begin{array}{ccc}
P_{m-2} & 0 & R_{m-2} \\
P_{m-1} & 0 & R_{m-1} \\
P_m & 1 & R_m \\
\end{array}\right|,\\
\nu_m&=\frac{(-1)^{m-1}}{q^{3m-3}}\left| \begin{array}{ccc}
P_{m-2} & Q_{m-2} & 0 \\
P_{m-1} & Q_{m-1} & 0 \\
P_m & Q_m & 1\\
\end{array}\right|.
\end{align*}
According to \eqref{A117-A118-A119F1} and \eqref{A117-A118-A119F2},
by setting
\begin{equation*}
\left\{\begin{array}{l} A_m=q^m(\lambda_m+q^2\lambda_{m+1}),\\
B_m=q^m(\mu_m+q^2\mu_{m+1}),\\
C_m=q^m(\nu_m+q^2\nu_{m+1}),
\end{array}\right. \quad \text{and} \quad
\left\{\begin{array}{l} E_m=\lambda_m+(q^2-q)\lambda_{m+1}-q^3\lambda_{m+2},\\
F_m=\mu_m+(q^2-q)\mu_{m+1}-q^3\mu_{m+2},\\
G_m=\nu_m+(q^2-q)\nu_{m+1}-q^3\nu_{m+2},
\end{array}\right.
\end{equation*}
we have
\begin{align*}
\sum_{n=0}^{\infty}\frac{(-q;q^2)_{n+1}q^{n^2+2mn}}{(q;q^2)_{2n+1}}&=
q^{-m}A_mP_{\infty}+q^{-m}B_mQ_{\infty}+q^{-m}C_mR_{\infty},\\
\sum_{n=0}^{\infty}\frac{(-q;q^2)_nq^{n^2+2mn}}{(q;q^2)_{2n}}&=
E_mP_{\infty}+F_mQ_{\infty}+G_mR_{\infty}.
\end{align*}
Since
\begin{equation*}
A_m=\frac{(-1)^{m-1}}{q^{2m-2}}\left| \begin{array}{ccc}
-1 & qQ_{m-2} & qR_{m-2} \\
0 & Q_{m-1} & R_{m-1} \\
1 & Q_m & R_m \\
\end{array}\right|,
\end{equation*}
by letting the last two columns in the determinants of $A_{m-1}$,
$A_{m-2}$, and $A_{m-3}$ be the same as those of $A_m$, we find a
linear equation
\begin{equation*}
A_m=A_{m-1}+(q^2+q^{2m-4})A_{m-2}-q^2A_{m-3}.
\end{equation*}
Using the initial conditions of $P_n$, $Q_n$, and $R_n$, we have
\begin{equation*}
A_0=1, \ A_1=0, \ A_2=0.
\end{equation*}
Proceeding in the same way, we get the recursions of $B_m$, $C_m$,
$E_m$, $F_m$, and $G_m$. Therefore, we obtain
\eqref{A117-A118-A119-1} and \eqref{A117-A118-A119-2}. \qed

The identity \eqref{A119} is a special case of
\eqref{A117-A118-A119-1}, and \eqref{A117} and \eqref{A118} are the
special cases of \eqref{A117-A118-A119-2}.

\begin{thm}\label{A21-Zimmer2.5-Bowman2.7Thm}We have
\begin{itemize}
\item[(1)]
\begin{align}\label{A21-Zimmer2.5-Bowman2.7-1}
\sum_{n=0}^{\infty}\frac{(-1)^n(q;q^2)_nq^{n^2+2mn}}{(-q;q^2)_{n+1}(q^4;q^4)_n}&=
\frac{(-1)^mq^{-m}(-q^2,-q^3,q^5;q^5)_{\infty}(q;q^2)_{\infty}}{(q^2;q^2)_{\infty}}A_m
\nonumber\\
&+\frac{(-1)^mq^{-m}(q^{10};q^{10})_{\infty}(q^{20};q^{20})_{\infty}}{(q;q^2)_{\infty}(q^5;q^{20})_{\infty}(q^4;q^4)_{\infty}}B_m
\nonumber\\
&+\frac{(-1)^mq^{-m}(-q,-q^4,q^5;q^5)_{\infty}(q;q^2)_{\infty}}{(q^2;q^2)_{\infty}}C_m,
\end{align}
where
\begin{align*}
&A_m=(1+q^{2m-4})A_{m-1}+(q^2+q^{2m-4})A_{m-2}-q^2A_{m-3},\qquad A_0=1, \ A_1=0, \ A_2=0,\\
&B_m=(1+q^{2m-4})B_{m-1}+(q^2+q^{2m-4})B_{m-2}-q^2B_{m-3},\qquad B_0=-q, \ B_1=-q, \ B_2=-q,\\
&C_m=(1+q^{2m-4})C_{m-1}+(q^2+q^{2m-4})C_{m-2}-q^2C_{m-3},\qquad
C_0=0, \ C_1=0, \ C_2=q.
\end{align*}
\item[(2)]
\begin{align}\label{A21-Zimmer2.5-Bowman2.7-2}
\sum_{n=0}^{\infty}\frac{(-1)^n(q;q^2)_nq^{n^2+2mn}}{(-q;q^2)_n(q^4;q^4)_n}&=
\frac{(-q^2,-q^3,q^5;q^5)_{\infty}(q;q^2)_{\infty}}{(q^2;q^2)_{\infty}}E_m
+\frac{(q^{10};q^{10})_{\infty}(q^{20};q^{20})_{\infty}}{(q;q^2)_{\infty}(q^5;q^{20})_{\infty}(q^4;q^4)_{\infty}}F_m
\nonumber\\
&+\frac{(-q,-q^4,q^5;q^5)_{\infty}(q;q^2)_{\infty}}{(q^2;q^2)_{\infty}}G_m,
\end{align}
where
\begin{align*}
&E_m=-(q+q^{2m-3})E_{m-1}+(1+q^{2m-4})E_{m-2}+qE_{m-3},\qquad E_0=1, \ E_1=0, \ E_2=1,\\
&F_m=-(q+q^{2m-3})F_{m-1}+(1+q^{2m-4})F_{m-2}+F_{m-3},\qquad F_0=0, \ F_1=0, \ F_2=2q,\\
&E_m=-(q+q^{2m-3})G_{m-1}+(1+q^{2m-4})G_{m-2}+qG_{m-3},\qquad G_0=0,
\ G_1=1, \ G_2=-q.
\end{align*}
\end{itemize}
\end{thm}
\pf The identity A.21 in Slater's list is stated as follows.

\textbf{Identity A.21 (Slater \cite{Slater}):}
\begin{equation}\label{A21}
\sum_{n=0}^{\infty}\frac{(-1)^n(q;q^2)_nq^{n^2}}{(-q;q^2)_n(q^4;q^4)_n}=
\frac{(-q^2,-q^3,q^5;q^5)_{\infty}(q;q^2)_{\infty}}{(q^2;q^2)_{\infty}}.
\end{equation}
\begin{align}\label{A21P}
&P_n=(1-q-q^2-q^{2n-1})P_{n-1}+(q+q^2-q^3+q^{2n-2})P_{n-2}+q^3P_{n-3},\nonumber
\\
&\qquad \qquad P_0=1, \ P_1=1-q, \ P_2=1-q+2q^2+q^4.
\end{align}

Recently, McLaughlin et al. and Bowman et al. found two new
Rogers-Ramanujan type identities in \cite{McLaughlin-Sills-Zimmer}
and \cite{Bowman-McLaughin-Sills}, respectively.

\textbf{An identity (McLaughlin et al. \cite[Eq.
(2.5)]{McLaughlin-Sills-Zimmer}):}
\begin{equation}\label{Zimmer2.5}
\sum_{n=0}^{\infty}\frac{(-1)^n(q;q^2)_nq^{n^2+2n}}{(-q;q^2)_{n+1}(q^4;q^4)_n}=
\frac{(q^{10};q^{10})_{\infty}(q^{20};q^{20})_{\infty}}{(q;q^2)_{\infty}(q^5;q^{20})_{\infty}(q^4;q^4)_{\infty}}.
\end{equation}

\textbf{An identity (Bowman et al. \cite[Thm.
2.7]{Bowman-McLaughin-Sills}):}
\begin{equation}\label{Bowman2.7}
\sum_{n=0}^{\infty}\frac{(-1)^n(q;q^2)_nq^{n^2+2n}}{(-q;q^2)_n(q^4;q^4)_n}=
\frac{(-q,-q^4,q^5;q^5)_{\infty}(q;q^2)_{\infty}}{(q^2;q^2)_{\infty}}.
\end{equation}
We can see that \eqref{Zimmer2.5} and \eqref{Bowman2.7} are partners
to \eqref{A21}. Therefore, we have
\begin{align}\label{Zimmer2.5Q}
&Q_n=(1-q-q^2-q^{2n-1})Q_{n-1}+(q+q^2-q^3+q^{2n-2})Q_{n-2}+q^3Q_{n-3},\nonumber
\\
&\qquad \qquad Q_0=0, \ Q_1=1, \ Q_2=1-q-q^3,
\end{align}
\begin{align}\label{Bowman2.7R}
&R_n=(1-q-q^2-q^{2n-1})R_{n-1}+(q+q^2-q^3+q^{2n-2})R_{n-2}+q^3R_{n-3},\nonumber
\\
&\qquad \qquad R_0=1, \ R_1=1, \ R_2=1-q^3,
\end{align}
where $P_n$, $Q_n$, and $R_n$ converge to the right hand sides of
\eqref{A21}, \eqref{Zimmer2.5}, and \eqref{Bowman2.7}, respectively.
The initial conditions for $Q_n$ and $R_n$ are obtained in the
following analysis.

Now we consider the following determinant:{\footnotesize
\begin{equation*}
F(z):=\left| \begin{array}{cccccc} 1-q-q^2-zq & q+q^2-q^3+zq^2 & q^3 & & & \cdots \\
-1 & 1-q-q^2-zq^3 & q+q^2-q^3+zq^4 & q^3 & & \cdots  \\
& -1 & 1-q-q^2-zq^5 & q+q^2-q^3+zq^6 & q^3 & \cdots \\
& & \ddots & \ddots & \ddots & \ddots\\
\end{array}
\right|.
\end{equation*}}
Expanding the determinant with respect to the first column, we get
\begin{equation*}
F(z)=(1-q-q^2-zq)F(zq^2)+(q+q^2-q^3+zq^2)F(zq^4)+q^3F(zq^6).
\end{equation*}
Setting
\begin{equation*}
F(z)=\sum_{n=0}^{\infty}a_nz^n,
\end{equation*}
we get, upon comparing coefficients,
\begin{equation*}
a_n=\frac{-(1-q^{2n-1})q^{2n-1}}{(1-q^{2n})(1+q^{2n+1})(1+q^{2n+2})}a_{n-1}=\cdots=
\frac{(-1)^n(q;q^2)_nq^{n^2}(1+q)(1+q^2)}{(-q;q^2)_{n+1}(q^4;q^4)_n(1+q^{2n+2})}a_{0}.
\end{equation*}
Since $a_0=\frac{1}{(1+q)(1+q^2)}$, we have
\begin{equation*}
F(z)=\sum_{n=0}^{\infty}\frac{(-1)^n(q;q^2)_nq^{n^2}}{(-q;q^2)_{n+1}(q^4;q^4)_n(1+q^{2n+2})}z^n.
\end{equation*}
Thus we get
\begin{equation}\label{A21-Zimmer2.5-Bowman2.7F1}
\sum_{n=0}^{\infty}\frac{(-1)^n(q;q^2)_nq^{n^2+2mn}}{(-q;q^2)_{n+1}(q^4;q^4)_n}=F(q^{2m})+q^2F(q^{2m+2}),
\end{equation}
\begin{equation}\label{A21-Zimmer2.5-Bowman2.7F2}
\sum_{n=0}^{\infty}\frac{(-1)^n(q;q^2)_nq^{n^2+2mn}}{(-q;q^2)_n(q^4;q^4)_n}=F(q^{2m})+(q^2+q)F(q^{2m+2})+q^3F(q^{2m+4}).
\end{equation}

On the other hand, $F(z)$ is the limit of the finite determinant
{\scriptsize
\begin{align*}\nonumber
D_n(z):= \left| \begin{array}{ccccc} 1-q-q^2-zq & q+q^2-q^3+zq^2 & q^3 & & \cdots \\
-1 & 1-q-q^2-zq^3 & q+q^2-q^3+zq^4 & q^3 & \cdots  \\
\vdots & \ddots & \ddots & \ddots & \ddots \\
& & -1 & 1-q-q^2-zq^{2n-3} & q+q^2-q^3+zq^{2n-2} \\
& & & -1 & 1-q-q^2-zq^{2n-1} \\
\end{array} \right|.
\end{align*}}
Expanding this determinant with respect to the last row, we get
\begin{align*}
&D_n(z)=(1-q-q^2-zq^{2n-1})D_{n-1}(z)+(q+q^2-q^3+zq^{2n-2})D_{n-2}(z)+q^3D_{n-3}(z),\\
&\qquad\qquad D_{-1}(z)=0,\ D_0(z)=1,\ D_1(z)=1-q-q^2-zq.
\end{align*}
Then we have{\footnotesize
\begin{equation*}
D_{n-m}(q^{2m})=(1-q-q^2-q^{2n-1})D_{n-m-1}(q^{2m})+(q+q^2-q^3+q^{2n-2})D_{n-m-2}(q^{2m})+q^3D_{n-m-3}(q^{2m}).
\end{equation*}}
Now we calculate the initial conditions of $Q_n$ and $R_n$ in
\eqref{Zimmer2.5Q} and \eqref{Bowman2.7R}. According to
\eqref{A21-Zimmer2.5-Bowman2.7F1} and
\eqref{A21-Zimmer2.5-Bowman2.7F2}, we have
\begin{align*}
Q_{\infty}&=F(q^2)+q^2F(q^4),\\
R_{\infty}&=F(q^2)+(q^2+q)F(q^4)+q^3F(q^6).
\end{align*}
Due to $\displaystyle\lim_{n \rightarrow
\infty}D_{n-m}(q^{2m})=F(q^{2m})$, we have
\begin{align*}
Q_n&=D_{n-1}(q^2)+q^2D_{n-2}(q^4),\\
R_n&=D_{n-1}(q^2)+(q^2+q)D_{n-2}(q^4)+q^3D_{n-3}(q^6).
\end{align*}
Therefore, we get
\begin{align*}
&Q_0=0, \ Q_1=1, \ Q_2=1-q-q^3;\\
&R_0=1, \ R_1=1, \ R_2=1-q^3.
\end{align*}

Since $\langle D_{n-m}(q^{2m})\rangle_n$, $\langle P_n\rangle_n$,
$\langle Q_n\rangle_n$, and $\langle R_n\rangle_n$ satisfy the
same recursion, we set
\begin{equation*}
D_{n-m}(q^{2m})=\lambda_m P_n+\mu_m Q_n+\nu_m R_n.
\end{equation*}
Using the initial conditions $D_{-1}(q^{2m})=0$, $D_{0}(q^{2m})=1$,
and $D_{1}(q^{2m})=1-q-q^2-q^{2m+1}$, we have
\begin{align*}
\lambda_m&=\frac{\left| \begin{array}{ccc}
0 & Q_{m-1} & R_{m-1} \\
1 & Q_m & R_m \\
1-q-q^2-q^{2m+1} & Q_{m+1} & R_{m+1} \\
\end{array}\right|}{\left|\begin{array}{ccc}
P_{m-1} & Q_{m-1} & R_{m-1} \\
P_m & Q_m & R_m \\
P_{m+1} & Q_{m+1} & R_{m+1} \\
\end{array}\right|},\\
\mu_m&=\frac{\left| \begin{array}{ccc}
P_{m-1} & 0 &  R_{m-1} \\
P_m & 1 &  R_m \\
P_{m+1} & 1-q-q^2-q^{2m+1} & R_{m+1} \\
\end{array}\right|}{\left|\begin{array}{ccc}
P_{m-1} & Q_{m-1} & R_{m-1} \\
P_m & Q_m & R_m \\
P_{m+1} & Q_{m+1} & R_{m+1} \\
\end{array}\right|},\\
\nu_m&=\frac{\left| \begin{array}{ccc}
P_{m-1} & Q_{m-1} & 0  \\
P_m & Q_m & 1   \\
P_{m+1} & Q_{m+1} & 1-q-q^2-q^{2m+1} \\
\end{array}\right|}{\left|\begin{array}{ccc}
P_{m-1} & Q_{m-1} & R_{m-1} \\
P_m & Q_m & R_m \\
P_{m+1} & Q_{m+1} & R_{m+1} \\
\end{array}\right|},
\end{align*}
where
\begin{equation*}
\left|\begin{array}{ccc}
P_{m-1} & Q_{m-1} & R_{m-1} \\
P_m & Q_m & R_m \\
P_{m+1} & Q_{m+1} & R_{m+1} \\
\end{array}\right|=-q^{3m-1}(1+q),
\end{equation*}
which can be proved by induction on $m$. Therefore, we have simpler
forms for $\lambda_m$, $\mu_m$, and $\nu_m$ as follows:
\begin{align*}
\lambda_m&=-\frac{\left| \begin{array}{ccc}
0 & Q_{m-2} & R_{m-2} \\
0 & Q_{m-1} & R_{m-1} \\
1 & Q_m & R_m \\
\end{array}\right|}{q^{3m-4}(1+q)},\\
\mu_m&=-\frac{\left| \begin{array}{ccc}
P_{m-2} & 0 & R_{m-2} \\
P_{m-1} & 0 & R_{m-1} \\
P_m & 1 & R_m \\
\end{array}\right|}{q^{3m-4}(1+q)},\\
\nu_m&=-\frac{\left| \begin{array}{ccc}
P_{m-2} & Q_{m-2} & 0 \\
P_{m-1} & Q_{m-1} & 0 \\
P_m & Q_m & 1\\
\end{array}\right|}{q^{3m-4}(1+q)}.
\end{align*}
According to \eqref{A117-A118-A119F1} and \eqref{A117-A118-A119F2},
by setting
\begin{equation*}
\left\{\begin{array}{l} A_m=(-1)^mq^m(\lambda_m+q^2\lambda_{m+1}),\\
B_m=(-1)^mq^m(\mu_m+q^2\mu_{m+1}),\\
C_m=(-1)^mq^m(\nu_m+q^2\nu_{m+1}),
\end{array}\right. \quad \text{and} \quad
\left\{\begin{array}{l} E_m=\lambda_m+(q^2+q)\lambda_{m+1}+q^3\lambda_{m+2},\\
F_m=\mu_m+(q^2+q)\mu_{m+1}+q^3\mu_{m+2},\\
G_m=\nu_m+(q^2+q)\nu_{m+1}+q^3\nu_{m+2},
\end{array}\right.
\end{equation*}
we have
\begin{align*}
\sum_{n=0}^{\infty}\frac{(-q;q^2)_{n+1}q^{n^2+2mn}}{(q;q^2)_{2n+1}}&=
(-1)^mq^{-m}A_mP_{\infty}+(-1)^mq^{-m}B_mQ_{\infty}+(-1)^mq^{-m}C_mR_{\infty},\\
\sum_{n=0}^{\infty}\frac{(-q;q^2)_nq^{n^2+2mn}}{(q;q^2)_{2n}}&=
E_mP_{\infty}+F_mQ_{\infty}+G_mR_{\infty}.
\end{align*}
Since
\begin{equation*}
A_m=\frac{(-1)^{m-1}}{q^{2m-3}(1+q)}\left| \begin{array}{ccc}
1 & qQ_{m-2} & qR_{m-2} \\
0 & Q_{m-1} & R_{m-1} \\
1 & Q_m & R_m \\
\end{array}\right|,
\end{equation*}
we can find a linear equation
\begin{equation*}
A_m=(1+q^{2m-4})A_{m-1}+(q^2+q^{2m-4})A_{m-2}-q^2A_{m-3}.
\end{equation*}
Using the initial conditions of $P_n$, $Q_n$, and $R_n$, we have
\begin{equation*}
A_0=1, \ A_1=0, \ A_2=0.
\end{equation*}
Proceeding in the same way, we get the recursions of $B_m$, $C_m$,
$E_m$, $F_m$, and $G_m$. Therefore, we obtain
\eqref{A21-Zimmer2.5-Bowman2.7-1} and
\eqref{A21-Zimmer2.5-Bowman2.7-2}. \qed

The identity \eqref{Zimmer2.5} is a special case of
\eqref{A21-Zimmer2.5-Bowman2.7-1}, and the identities \eqref{A21}
and \eqref{Bowman2.7} are the special cases of
\eqref{A21-Zimmer2.5-Bowman2.7-2}.

%================================================================

\end{document}